\documentclass{amsart}

\usepackage{amsmath}
\usepackage{amssymb}
\usepackage{amsthm}
\usepackage{xypic}
\usepackage{algorithmic}
\usepackage{algorithm}
\usepackage[pdftex]{graphicx}

\theoremstyle{plain}
\newtheorem{theorem}{Theorem}[section]
\newtheorem*{theorem*}{Theorem}
\newtheorem*{corollary*}{Corollary}
\newtheorem{lemma}[theorem]{Lemma}
\newtheorem{proposition}[theorem]{Proposition}
\newtheorem{corollary}[theorem]{Corollary}

\theoremstyle{definition}
\newtheorem{definition}[theorem]{Definition}

\numberwithin{equation}{section}
\numberwithin{table}{section}

\def\Z{\mathbb{Z}}
\def\F{\mathbb{F}}
\def\Q{\mathbb{Q}}
\def\C{\mathbb{C}}

\def\OK{\mathcal{O}_K}

\def\qed{\raisebox{-.25ex}{\scalebox{.786}[1.272]{$\blacksquare$}}}
\def\divides{\,|\,}

\DeclareMathOperator{\Gal}{Gal}

\makeatletter
\def\pmod#1{\allowbreak\mkern10mu({\operator@font mod}\,\,#1)}
\def\pod#1{\allowbreak\mkern10mu(#1)}
\makeatother

\newcommand{\sumstar}{\sideset{}{^\star}\sum}
\newcommand{\psimodq}{\substack{\psi\\\text{mod }q }}

\title{Norm-Euclidean Galois Fields}
\author{Kevin J. McGown}
\date{}
\address{Department of Mathematics, University of California, San Diego, 9500 Gilman Drive, La Jolla, CA 92093}
\curraddr{Department of Mathematics, Oregon State University, 368 Kidder Hall, Corvallis, OR 97331}
\email{mcgownk@math.oregonstate.edu}
\keywords{norm-Euclidean, Galois fields, cubic fields, character sums}
\subjclass[2010]{Primary 11A05, 11R04, 11Y40; Secondary 11R16, 11R80, 11L40, 11R32}

\begin{document}

\maketitle

\begin{abstract}
Let K be a Galois number field of prime degree $\ell$.
Heilbronn showed that for a given $\ell$ there are only finitely many such fields that are norm-Euclidean.
In the case of $\ell=2$ all such norm-Euclidean fields have been identified, but for $\ell\neq 2$,
little else is known.  We give the first upper bounds on the discriminants of such fields when $\ell>2$.
Our methods lead to a simple algorithm which allows one to generate a list of candidate
norm-Euclidean fields up to a given discriminant, and we provide some computational results.
\end{abstract}



\section{Introduction}\label{S:intro}

Let $K$ be a number field with ring of integers $\OK$, and denote by $N=N_{K/\Q}$ the
absolute norm map.  For brevity, we will sometimes use the term field to mean a number field.
We call a number field $K$ norm-Euclidean if for every $\alpha,\beta\in\OK$, $\beta\neq 0$, there exists
$\gamma\in\OK$ such that $|N(\alpha-\gamma\beta)|<|N(\beta)|$.
In the quadratic setting, it is known that there are only finitely many norm-Euclidean fields and they have been identified;
namely, a number field of the form $K=\Q(\sqrt{d})$ with $d$ squarefree is norm-Euclidean if and only if
$$
  d=-1,-2,-3,-7,-11,2, 3, 5, 6, 7, 11, 13, 17, 19, 21, 29, 33, 37, 41, 57, 73
  \,.
$$
 In his third and final paper on the Euclidean algorithm (see~\cite{heilbronn:cyclic}),
Heilbronn proves a finiteness result for various classes of
cyclic fields.  For us, the most important part of Heilbronn's result states:
\begin{theorem}[Heilbronn, 1951]\label{T:Heilbronn}
Given a prime $\ell$, there are only finitely many norm-Euclidean Galois fields
of degree $\ell$.
\end{theorem}
However, Heilbronn's result on cyclic fields
does not give an upper bound on the discriminant, even in the cubic case.
The case of Galois cubic fields is especially interesting, as we have the following (see~\cite{godwin:1965, smith:1969, godwin.smith:1993}):
\begin{theorem}[Godwin \& Smith, 1993]\label{T:godwin.smith}
The norm-Euclidean Galois cubic fields
with discriminant \mbox{$|\Delta|<10^8$} are exactly those with
$$
    \Delta=7^2,9^2,13^2,19^2,31^2,37^2,43^2,61^2,67^2,103^2,109^2,127^2,157^2
    \,.
$$  
\end{theorem}
Lemmermeyer has further verified that this list constitutes all fields with $|\Delta|<2.5\cdot 10^{11}$ (see~\cite{lemmermeyer:euclidean}).
We prove the following result, which gives an upper bound on the discriminant for the fields
considered in Theorem~\ref{T:Heilbronn}.
\begin{theorem}\label{T:DB}
  Let $\ell$ be an odd prime.  There exists
  computable constant $C_\ell$ such that if $K$ is
  a Galois number field of odd prime degree $\ell$, conductor $f$,
  and discriminant $\Delta$, which is norm-Euclidean,
  then $f<C_\ell$ and $0<\Delta<C_\ell^{\ell-1}$.
\begin{table}[H]
\centering
\begin{tabular}{c c}
\begin{tabular}{| l | l |}
\hline
$\ell$ & $\;C_\ell$ \\
\hline
$3$ & $10^{70}$\\
$5$ & $10^{78}$\\
$7$ & $10^{82}$\\
$11$ & $10^{88}$\\
$13$ & $10^{89}$\\
$17$ & $10^{92}$\\
$19$ & $10^{94}$\\
$23$ & $10^{96}$\\
\hline
\end{tabular}
\qquad
\begin{tabular}{| l | l |}
\hline
$\ell$ & $\;C_\ell$ \\
\hline
$29$ & $10^{98}$\\
$31$ & $10^{99}$\\
$37$ & $10^{101}$\\
$41$ & $10^{102}$\\
$43$ & $10^{102}$\\
$47$ & $10^{103}$\\
$53$ & $10^{104}$\\
$59$ & $10^{105}$\\
\hline
\end{tabular}
\qquad
\begin{tabular}{| l | l |}
\hline
$\ell$ & $\;C_\ell$ \\
\hline
$61$ & $10^{106}$\\
$67$ & $10^{107}$\\
$71$ & $10^{107}$\\
$73$ & $10^{108}$\\
$79$ & $10^{108}$\\
$83$ & $10^{109}$\\
$89$ & $10^{109}$\\
$97$ & $10^{110}$\\
\hline
\end{tabular}
\end{tabular}
\caption{Values of $C_\ell$ for primes $\ell<100$\label{Table:Cl}}
\end{table}
\end{theorem}

Although the results of the previous theorem represent a significant step forward,
the magnitude of the constants leaves something to be desired, especially
if one is interested in determining all such fields, for any fixed $\ell$.
As is frequently the case in
estimates of number theoretic quantities,
under the Generalized Riemann Hypothesis (GRH) one should be able to obtain much sharper results.
This is the subject of a forthcoming paper (see~\cite{mcgown:grh}).

In order to prove Theorem~\ref{T:DB}, we derive explicit inequalities which guarantee the failure of the norm-Euclidean property.
Our inequalities (see Theorem~\ref{T:conditions}) involve the existence of small integers satisfying certain splitting and congruence conditions.
This also leads to an algorithm (see Algorithm~\ref{Alg:1}) for tabulating a list of candidate norm-Euclidean Galois fields (of prime degree $\ell$) up to a given discriminant.
We have implemented this algorithm in the mathematics software SAGE, thereby obtaining:
\begin{theorem}\label{T:comp}
The following table contains all possible norm-Euclidean Galois number fields of prime degree $\ell$ and conductor $f$ with $3\leq\ell\leq 30$ and $f\leq 10^4$.
(Of course, some of these fields may not be norm-Euclidean.)
\begin{table}[H]
\label{Table:candidate}
\centering
\begin{tabular}{| l | l |}
  \hline
  $\ell$ & $f\leq 10^4$ \\
  \hline
  $3$ & $7$, $9$, $13$, $19$, $31$, $37$, $43$, $61$, $67$, $73$, $103$, $109$, $127$, $157$,\\
  &  $277$, $439$, $643$, $997$, $1597$ \\ 
  $5$ & $11$, $25$, $31$, $41$, $61$, $71$, $151$, $311$, $431$ \\ 
  $7$ & $29$, $43$, $49$, $127$, $239$, $673$, $701$, $911$ \\
  $11$ & $23$, $67$, $89$, $121$, $331$, $353$, $419$, $617$ \\
  $13$ & $53$, $79$, $131$, $157$, $169$, $313$, $443$, $521$, $937$ \\
  $17$ & $137$, $289$, $443$, $1259$, $2687$ \\
  $19$ &  $191$, $229$, $361$, $1103$\\
  $23$ & $47$, $139$, $277$, $461$, $529$, $599$, $691$, $967$, $1013$, $1289$\\
  $29$ & $59$, $233$, $523$, $841$, $929$, $2843$, $3191$\\
  \hline
\end{tabular}
\caption{Candidate norm-Euclidean fields of small degree}
\end{table}
\end{theorem}
Notice that when $\ell=3$, we cover all possible $|\Delta|<10^8$
(as $\Delta=f^2$ in this case) and that
our results are consistent with Theorem~\ref{T:godwin.smith}.
We have opted not to remove those fields which are known not to be norm-Euclidean
from our table, but rather to give the exact output produced by our algorithm.
In the case of $\ell=3$, we know that exactly $13$ of the fields listed are norm-Euclidean,
but not too much seems to be known about the remaining fields in the table.
It would be interesting to study these fields using other methods,
possibly on a case-by-case basis if necessary,
to decide which among them are norm-Euclidean.\footnote{Of course, to begin with, one could determine which have class number one.}
No effort has been made in this direction by the author, but this may be the subject of a future investigation.

%
%

For the cubic case, we have implemented an efficient version
of the algorithm in the programming language C, which takes advantage of the cubic reciprocity law.
After 91 hours of computation on an iMac, we have verified that Godwin and Smith's list
constitutes all norm-Euclidean Galois cubic fields with $|\Delta|<10^{20}$.
Combining this with Theorem~\ref{T:DB} leads to
the following result
which represents the current state of knowledge for norm-Euclidean Galois cubic fields:
\begin{theorem}\label{T:comp.cubic}
The Galois cubic fields with
$$
    \Delta=7^2,9^2,13^2,19^2,31^2,37^2,43^2,61^2,67^2,103^2,109^2,127^2,157^2
    \,.
$$  
are norm-Euclidean,
and any remaining
norm-Euclidean Galois cubic field
must have discriminant
$\Delta=f^2$ with $f\equiv 1\pmod{3}$
where $f$
is a prime in the interval $(10^{10},\,10^{70})$.
\end{theorem}
Finally, we mention that under the GRH we can significantly improve on the above result (see~\cite{mcgown:grh}).


\section{Preliminaries}\label{S:prelim}

%
%

\subsection{Class field theory}\label{S:prelim.classfieldtheory}

We review some well-known facts regarding class field theory over $\Q$ which will be useful in the sequel;
one reference (among many possibilities) is~\cite{garbanati:cft}.
The famous Kronecker--Weber Theorem states that every abelian extension $K/\Q$ is contained in a cyclotomic extension,
and the conductor of $K$ is defined to be the smallest $f\in\Z^+$ such that $K\subseteq\Q(\zeta_f)$.  

Suppose $K/\Q$ is abelian and $K\subseteq\Q(\zeta_m)$.
We associate a character group $X_K$ to $K$ in the following manner.
Via Galois theory, we can identify $K$ with a subgroup $H$ of $(\Z/m\Z)^\star$,
and we define $X_K$ to be the subgroup of Dirichlet characters modulo $m$
that are trivial on $H$.
Different choices of $m$ lead to isomorphic
character groups $X_K$ in a natural way and one has $X_K\simeq\Gal(K/\Q)$.
Moreover,
the map $K\mapsto X_K$
gives a one-to-one correspondence between subfields of $\Q(\zeta_m)$
and subgroups of the character group of $(\Z/m\Z)^\star$.
Perhaps the most important property of this correspondence is that a rational prime $p$ splits in $K$
if and only if $\chi(p)=1$ for all $\chi\in X_K$.


In the case of interest to us,
$K/\Q$ is a cyclic number field of degree~$\ell$.
Suppose $K$ has conductor $f$, and view $X_K$ as a subgroup of
the group of Dirichlet characters modulo $f$.
In this case $X_K$ is cyclic and any generator 
is a primitive character modulo $f$ of order $\ell$.
Hence we have the following one-to-one correspondence:
\begin{equation}
  \nonumber
  \left\{
  \begin{array}{cc}
  \text{Cyclic extensions}
   \\
   \text{ $K/\Q$ of conductor $f$}
   \\
   \text{ and degree $\ell$}
   \end{array}
  \right\}
  \longleftrightarrow
  \left\{
  \begin{array}{cc}
  \text{Primitive Dirichlet characters}\\
  \text{$\chi:\Z\to\C$ of modulus $f$}
  \\
  \text{and order $\ell$}
  \end{array}
  \right\}
  \scalebox{1.5}{$/\sim$}
\end{equation}
The equivalence is given by the natural action of 
$\Gal(\overline{\Q}/\Q)$;
namely, $\chi\sim\psi$ if $\sigma\circ\chi=\psi$ for some \mbox{$\sigma\in\Gal(\overline{\Q}/\Q)$}.
If we are considering fields of degree $\ell$ it suffices to only consider those $\sigma\in\Gal(\Q(\zeta_\ell)/\Q)$,
as Dirichlet characters of order~$\ell$ take values in $\Q(\zeta_\ell)$,
and so this equivalence amounts to a choice of a primitive $\ell$-th root of unity among the $\phi(\ell)$ possibilities.
Moreover, this correspondence is such that a rational prime $p$ splits in $K$ if and only if $\chi(p)=1$.



\subsection{Number fields with class number one}\label{S:prelim.classnumberone}

In trying to locate Euclidean fields we may
restrict our attention to the case where $K$ has class number one,
so it is useful to see what
extra conditions this places on our fields.
We denote by $h_K$ and $h_K^+$ the class number and narrow class number respectively.
The next two lemmas are well-known, but we provide the proofs as they
can be difficult to find in the literature.

The result contained in the following lemma is
perhaps most elegantly demonstrated via genus theory (see~\cite{ishida:book}),
and so we first recall some definitions.
The genus field of an abelian number field $K$, denoted by $K^\text{g}$, is the largest absolutely abelian extension of $K$
that is unramified at all finite primes, and the genus number of $K$ is defined by $g_K:=[K^\text{g}:K]$.
It is well known (and follows immediately from class field theory)  that $g_K$
divides $h_K^+$.


\begin{lemma}\label{L:primepower}
  Suppose $K/\Q$ is cyclic with odd prime degree $\ell$ and discriminant $\Delta$.  If $t$ distinct rational primes
  divide $\Delta$, then $\ell^{t-1}$ divides $h_K$.
\end{lemma}

\noindent\textbf{Proof.}
By Theorem 5 of \cite{ishida:book} we have
$g_K=\ell^{t-1}$.
We know $g_K$ divides $h_K^+$; moreover, since $\ell$ is an odd prime and
$h_K^+$ differs from $h_K$ only by a power of $2$, we conclude that $g_K$ divides $h_K$ as well.
\qed




\begin{lemma}\label{L:disc}
Suppose $K/\Q$ is cyclic with odd prime degree $\ell$, conductor $f$, and discriminant $\Delta$.
Further, suppose that $K$ has class number one.
In this case, one has
$\Delta=f^{\ell-1}$.  Moreover:
\begin{enumerate}
\item
If $\gcd(f,\ell)=1$, then $f$ is a prime with $f\equiv1\pmod\ell$.
\item
If $\gcd(f,\ell)>1$, then $f=\ell^2$.
\end{enumerate}
\end{lemma}

\noindent\textbf{Proof.}
Since $K$ has class number one, Lemma~\ref{L:primepower} allows us to conclude
that $|\Delta|$ is a prime power.
Also, note that $K$ is totally real (as it is Galois of odd degree)
and hence $\Delta>0$.  Since $K$ is cyclic of prime degree,
the conductor--discriminant formula allows us 
to conclude that $\Delta=f^{\ell-1}$.  It remains to determine $f$.

Since $f$ is the conductor of $K$, we have the following inclusion of fields:
$\Q\subset K\subset\Q(\zeta_f)$.
Since $\Delta$ is only divisible by one prime, $f$ must also be divisible by a single prime;
say, $f=p^k$ for some prime $p$ and $k\in\Z^+$.
Thus $[\Q(\zeta_f):\Q]=\phi(f)=(p-1)p^{k-1}$.
From the inclusion of fields, we see
$[K:\Q]$ divides $[\Q(\zeta_f):\Q]$;
that is, $\ell$ divides $(p-1)p^{k-1}$.

At this point, we break the proof into cases.
First, suppose that $\gcd(f,\ell)=1$, so that $\ell\neq p$.
In this case, we must have $\ell$ divides $p-1$; that is,
$p\equiv 1\pmod\ell$.  This implies, via Galois theory, that there exists a cyclic field
$K'$ as depicted below:
$$
  \xymatrixrowsep{1pc}
  
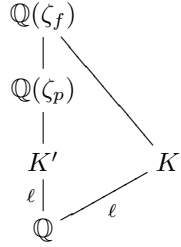
\xymatrix{
  \Q(\zeta_f)\ar@{-}[d]\ar@{-}[ddr]
  \\
  \Q(\zeta_p)\ar@{-}[d]
  \\
  K'\ar@{-}[d]_\ell & K\ar@{-}[dl]^\ell
  \\
  \Q
  }
$$
As there cannot be two cyclic fields of degree $\ell$ contained in $\Q(\zeta_f)$, we must have
that $K=K'$ and hence the conductor of $K$ satisfies $f=p$.
This proves the result in the case where $\gcd(f,\ell)=1$.

Turning to the case of $\gcd(f,\ell)>1$, we have $p=\ell$.
We argue as before, but this time we have $\ell$ divides $(\ell-1)\ell^{k-1}$;
this implies $k\geq 2$.  But we observe, as before, that
$\Q(\zeta_{\ell^2})$ already contains a cyclic field of order $\ell$, and hence
$f=\ell^2$,
which completes the proof.
\qed

\subsection{Heilbronn's criterion}\label{S:prelim.heilbronn}

Let $K$ be a Galois number field of odd prime degree $\ell$ and conductor $f$.
In order to state Heilbronn's criterion, we distinguish two subsets of the rational integers:
the norms,
$$
  \mathcal{N}:=N_{K/\Q}(\mathcal{O}_K)=\{n\in\Z\mid N_{K/\Q}(\alpha)=n \text{ for some $\alpha\in \OK$}\}
  \,,
$$
and the $\ell$-th power residues modulo $f$,
$$
  \mathcal{P}:=\{n\in\Z\mid x^\ell\equiv n\pmod f\text{ is soluble}\}
  \,.
$$
Although not stated in this way, Heilbronn proves the following~\cite{heilbronn:cyclic}:
\begin{lemma}[Heilbronn's Criterion]\label{L:NotEuclid}
Let $K$ be a Galois number field of odd prime degree $\ell$ and conductor $f$, with $(f,\ell)=1$.
If one can write $f=a+b$ with $a,b>0$, where $a,b\notin\mathcal{N}$ and $a\in\mathcal{P}$,
then $K$ is not norm-Euclidean.
\end{lemma}
This simple yet ingenious observation, which has its roots in a paper of
Erd\"os and Ko on quadratic fields~\cite{erdos.ko},
turns the problem into one of additive number theory.
For the sake of completeness, we provide the argument.

\vspace{1ex}
\noindent\textbf{Proof.}
Let $K$ be as in the hypothesis, and moreover, assume that $K$ is norm-Euclidean.
Suppose $f=a+b$ with $a,b>0$ where $a,b\notin\mathcal{N}$ and $a\in\mathcal{P}$.
We seek a contradiction.

Since $K$ is norm-Euclidean, it has class number one.
It follows from Lemma~\ref{L:disc} that $f$ is a prime,
and since $K$ has prime degree
we know that $f$ is totally ramified in $K$.
We factor $f=u\pi^\ell$ in $K$ where $\pi$ is a first degree prime and $u$ is a unit.
Fix an arbitrary $n\in\Z^+$.  There exists $\alpha\in \OK$ such that $n\equiv\alpha\pmod{\pi}$
with $|N(\alpha)|<|N(\pi)|=f$.
Conjugation gives $n\equiv\alpha^\sigma\pmod{\pi}$ for all embeddings $\sigma:K\to\C$,
and hence 
  $n^\ell\equiv N(\alpha)\pmod{f}$.
  Now we choose $n$ so that $a\equiv n^\ell\pmod{f}$ and we have
  $a\equiv N(\alpha)\pmod{f}$.
  Since $|N(\alpha)|<f$, we have either $N(\alpha)=a$ or $N(\alpha)=a-f=-b$.
  Thus $a$ or $-b$ lies in $\mathcal{N}$, a contradiction!
\qed
\vspace{1ex}


\section{Conditions for the Failure of the Norm-Euclidean Property}\label{S:conditions}

Building on the work of Heilbronn,
we prove the following theorem, which gives various conditions under which $K$ fails be norm-Euclidean.

\begin{theorem}\label{T:conditions}
  Let $K$ be a Galois number field of odd prime degree $\ell$ and conductor $f$ with $(f,\ell)=1$,
  and
  let $\chi$ be a primitive Dirichlet character modulo $f$ of order $\ell$.
  Denote by $q_1<q_ 2$ the two smallest rational primes that are inert in $K$.
  Suppose that there exists $r\in\Z^+$ with
  $$
    (r,q_1 q_2)=1,
    \quad
    \chi(r)=\chi(q_2)^{-1},
  $$
  such that any of the following conditions hold:
  \begin{enumerate}
    \item
    \qquad
    $r q_2 k\not\equiv f \pmod{q_1^2}$, \quad $k=1,\dots,q_1-1$,
    
    \quad\hspace{-1ex}
    $(q_1-1)(q_2 r-1)\leq f$
    
  \item
  \qquad
  $q_1\neq 2,3$,
  \quad
  $3q_1 q_2 r\log q_1< f $

  \item
  \qquad
  $q_1\neq 2,3,7$,
  \quad
  $2.1\,q_1 q_2 r\log q_1< f$
  
  \item
  \qquad
  $q_1=2,\,q_2\neq 3$,
  \quad
  $3q_2 r< f$
  
  \item
  \qquad
  $q_1=3,\,q_2\neq 5$,
  \quad
  $5q_2 r< f$
  \end{enumerate}
   Then $K$ is not norm-Euclidean.
\end{theorem}

The first condition in the above theorem places no restrictions on $q_1$ or $q_2$ but requires congruence conditions which hold ``most of the time'',
but can be rather awkward to verify.  The remaining conditions resulted from an effort to remove these congruences.

As in the statement of the above theorem, we will assume throughout this section
that $K$ is a Galois number field of odd prime degree $\ell$ and conductor $f$, with $(f,\ell)=1$,
and we will denote by $q_1<q_ 2$ the two smallest rational primes that are inert in $K$.
It suffices to assume that $K$ has class number one (otherwise it is immediate that $K$ is not norm-Euclidean), and we will do so.
Now Lemma~\ref{L:disc} tells us that the discriminant of $K$ satisfies $\Delta=f^{\ell-1}$,
where $f$ is a prime satisfying $f\equiv 1\pmod \ell$.
In light of Lemma~\ref{L:NotEuclid}, the following subset of $\Z^+$ will play a crucial role:
\begin{definition}
Let $\mathcal{S}$ denote the subset of positive integers less than $f$ which consists of
$\ell$-th power residues that are not norms.
In the notation of \S\ref{S:prelim.heilbronn}, \\
\mbox{$\mathcal{S}:=\mathcal{P}\cap\mathcal{N}^C\cap(0,f)$}.
\end{definition}
The following simple lemma characterizes $\mathcal{S}$ in terms
of~$\chi$, and will be used without comment in the arguments that follow.
\begin{lemma}
$$
       \mathcal{S}=\{n\in\Z\cap(0,f) \mid n=bc,\, (b,c)=1,\, \chi(b)\neq 1,\, \chi(bc)=1\}
$$
\end{lemma}

\noindent\textbf{Proof.}
Suppose $n\in\Z$ with $0<n<f$.
One knows that $n\in\mathcal{P}$
if and only if $\chi(n)=1$,
and that $n\notin\mathcal{N}$
if and only if one can write $n=bc$ with $(b,c)=1$ and $\chi(b)\neq 1$.
The result follows.
\qed

\begin{lemma}\label{L:preprop}
If there exists $s\in\mathcal{S}$ such that $(q_1,s)=1$ and $(q_1-1)(s-1)\leq f$, then
we can write
$f=us+v q_1$
with $0<u<q_1$ and $v>0$.
If $(q_1,v)=1$ in this expression, then $K$ is not norm-Euclidean.
\end{lemma}

\noindent\textbf{Proof.}
By a well-known theorem in elementary number theory,
the facts
$(q_1,s)=1$ and $(q_1-1)(s-1)\leq f$ imply
that there exists $u,v\in\Z_{\geq 0}$ such that
$f=us+v q_1$.
However, since $f$ is a prime not equal to $q_1$ and $s$ is composite,
we must have $u,v>0$, lest we arrive at a contradiction.
Without loss of generality, we can assume $u<q_1$.
Indeed, we just subtract multiples of $q_1$ from $u$ and add them to $v$ as necessary,
and the resulting $u$ and $v$ will remain positive for the same reason as before.
Since $u<q_1$, we have $\chi(p)=1$ for every prime $p$ dividing $u$,
and it follows that $us\in\mathcal{S}$.
If it were the case that $(q_1,v)=1$, then
we would have $v q_1\notin\mathcal{N}$ since $q_1\notin\mathcal{N}$;
in this case Lemma~\ref{L:NotEuclid} implies that $K$ is not norm-Euclidean.
\qed

\begin{proposition}\label{P:1}
If there exists $s\in\mathcal{S}$ such that $(s,q_1)=1$,
$$
  sk\not\equiv f\pmod{q_1^2}
  \,,\quad
  k=1,\dots,q_1-1
  \,,
$$
$$
  (q_1-1)(s-1)\leq f
  \,,
$$
then $K$ is not norm-Euclidean.
\end{proposition}

\noindent\textbf{Proof.}
By Lemma~\ref{L:preprop} we can write
$f=us+v q_1$ with $0<u<q_1$, $v>0$ and we may assume $q_1\divides v$.
This implies
$f\equiv us\pmod{q_1^2}$, a contradiction.
\qed


\vspace{1ex}
When $q_1\neq 2,3$, we can eliminate the congruence condition of Proposition~\ref{P:1}, but for a small price.

\begin{proposition}\label{P:2}
  Fix $q_1\neq 2,3$.  Suppose there exists a constant $1\leq B\leq 3$ such that
  for all $u\in(0,q_1)$ there exists a prime $p_0<B\log q_1$ with $(p_0,u)=1$.
  If there exists $s\in\mathcal{S}$ such that $(s,q_1)=1$ and
  $$
    Bq_1s\log q_1 \leq f
    \,,
  $$
  then $K$ is not norm-Euclidean.
\end{proposition}

\noindent\textbf{Proof.}
By Lemma~\ref{L:preprop} we can write
$f=us+vq_1$ with $0<u<q_1$, $v>0$
and we may assume
$q_1\divides v$.
By our hypothesis, there exists a prime such that $(p_0,u)=1$ and $p_0<B\log q_1$ for some $B\in[1,3]$.
In particular, we have $p_0<q_1$ since $3\log q_1<q_1$ for $q_1\geq 5$.
Let $n$ denote the smallest positive solution to the congruence
$$
  u+nq_1\equiv 0\pmod{p_0}
  \,,
 $$
so that $0<n<p_0$.
We claim that the expression
\begin{equation}\label{E:desired}
  f=(u+nq_1)s+(v-ns)q_1
\end{equation}
is of the desired form (to which Lemma~\ref{L:NotEuclid} applies).  First we note that
$$
  u+n q_1<q_1+(p_0-1)q_1=p_0 q_1
  \,.
$$  
To see that both terms in (\ref{E:desired}) are positive we
observe
$$
  (u+n q_1)s<p_0 q_1 s<B q_1 s \log q_1 \leq f
  \,.
$$
Notice that every prime $p$ dividing $u+n q_1$ is less than $q_1$,
which says $(u+n q_1)s\in\mathcal{S}$, as before.
If it were the case that $q_1 | v-ns$, then we would have
$q_1|s$, a contradiction; hence $(q_1, v-n s)=1$.
Now Lemma~\ref{L:NotEuclid} gives the result.
\qed

\vspace{1ex}
Motivated by the previous proposition, we introduce the following lemma which
gives the existence of the constant $B$.

\begin{lemma}\label{L:prime}
Suppose $q$ is prime and $0<u<q$.
If $q\neq 2,3$, then there exists a prime $p_0<3\log q$ such that $(p_0,u)=1$.
If $q\neq 2,3,7$, then there exists a prime $p_0<2.1\log q$ such that $(p_0,u)=1$.
\end{lemma}

\noindent\textbf{Proof.}
To show there exists a prime $p_0\leq x$ with $(p_0,u)=1$
it suffices to show
$$
  \sum_{p\leq x}\log p>\log u
  \,,
$$
as this implies the desired result.
For any $x\geq 5$ we have the inequality
\begin{equation}\label{E:RS}
  \sum_{p\leq x} \log p > \frac{x}{2.1}
  \,,
\end{equation}
which is easily deduced
from Corollary 3.16 of \cite{rosser.schoenfeld}
with a small amount of computation.\footnote
{
In fact, one can demonstrate this using the elementary methods
given in Ch. XXII of~\cite{hardy.wright}
together with an explicit version of Stirling's formula
if one is willing to do a little more computation.
}
Using this fact together with the hypothesis that $u<q$, one sees that
it suffices to show
\begin{equation}\label{E:crit}
  \log q\leq \frac{x}{2.1} 
  \,.
\end{equation}
This condition clearly holds when we set $x=2.1\log q$.
When $q\geq 11$, we have
$x\geq 2.1\log 11>5$, and the proof is complete.
The cases of $q=5,7$ are done by direct inspection.
\qed


\begin{proposition}\label{P:3}
  Suppose $q_1=2$, $q_2\neq 3$.
  If there exists $s\in\mathcal{S}$ such that
  \mbox{$(q_1,s)=1$}
  and $3s< f$,
  then $K$ is not norm-Euclidean.
\end{proposition}


\noindent\textbf{Proof.}
By Lemma~\ref{L:preprop} we may assume 
$f=s+2v$ with $2\divides v$.  In this case, we write
$f=3s+2(v-s)$.
If it were the case that $2 \divides (v-s)$, then we would have $2\divides s$, a contradiction.
Also observe that $\chi(3)=1$ and hence $3s\in\mathcal{S}$.
Finally, notice that
$3s<f$, which implies $v-s>0$.
\qed

\begin{proposition}\label{P:4}
  Suppose $q_1=3$, $q_2\neq 5$.
  If there exists $s\in\mathcal{S}$ such that
  \mbox{$(q_1,s)=1$}
  and
  $5s< f$,
  then $K$ is not norm-Euclidean.
\end{proposition}

\noindent\textbf{Proof.}
By Lemma~\ref{L:preprop} we may assume
$f=us+3v$ with $0<u<3$, $v>0$, and
$3\divides v$.
We treat separately the cases of $u=1$ and $u=2$.
If $u=1$, we have $f=s+3v$, which we rewrite as
$f=4s+3(v-s)$.
Proceeding as before we find this expression is of the desired form (since $\chi(2)=1$), provided $4s\leq f$.
If $u=2$, we have $f=2s+3v$, which we rewrite as
$f=5s+3(v-s)$,
which is of the desired form (since $\chi(5)=1$), provided $5s< f$.
\qed

\vspace{1ex}
Now we are ready:

\vspace{1ex}
\noindent\textbf{Proof of the theorem.}
If condition~(1) holds, we apply Proposition~\ref{P:1} with $s=q_2 r$.
If either of conditions~(2) or (3) hold, then we apply Proposition~\ref{P:2} with $s=q_2 r$
and invoke Lemma~\ref{L:prime}.
If conditions~(4) or (5) hold, we apply Propositions~\ref{P:3} or \ref{P:4} respectively.
\qed


\section{Discriminant Bounds}\label{S:DB}

\subsection{Some special cases}\label{S:DB.1}

The goal in \S\ref{S:DB.1} is to prove the following proposition which treats two very special cases.
The purpose of this is two-fold:  This will serve as an illustration of the type of inequalities we seek;
and, this will allow us to rid ourselves of these two cases which are particularly troublesome.

\begin{proposition}\label{P:special}
Let $K$ be a Galois number field of odd prime degree $\ell$ and conductor $f$.
Denote by $q_1<q_ 2$ the two smallest rational primes that are inert in $K$.
Suppose either of the following conditions hold:
\begin{enumerate}
  \item
    \qquad
    $q_1=2,\, q_2=3$,
    
    \quad\hspace{-1ex}
    $72(\ell-1)f^{1/2}\log 4f+35\leq f$
  \item
    \qquad
      $q_1=3,\, q_2=5$,
      
    \quad\hspace{-1ex}
      $507(\ell-1)f^{1/2}\log 9f+448\leq f$
\end{enumerate}
Then $K$ is not norm-Euclidean.
\end{proposition}

Notice that the above inequalities are completely explicit, they involve only $\ell$ and $f$, and for fixed $\ell$
they clearly hold beyond some easily computed value of $f$.
Ultimately, we will derive an analogous result which holds regardless of the values of $q_1$ and $q_2$
(see Theorem~\ref{T:DB.inequality}).
The following corollary, whose proof is immediate, 
is an example of the type of discriminant bound we can obtain from Proposition~\ref{P:special}.

\begin{corollary}
  Suppose $K$ is a norm-Euclidean Galois cubic field of conductor $f$ and discriminant $\Delta$.
  If the primes $2$ and $3$ are inert in $K$, then $f<10^7$ and $0<\Delta<10^{14}$.
\end{corollary}

First we prove a lemma about Dirichlet characters.

\begin{lemma}\label{L:DC1}
Suppose $\chi$ is a Dirichlet character
modulo $m$ of order $\ell$.
Fix an $\ell$-th root of unity $\zeta$.
Let ($\star$) be any property of integers.
Suppose there are no integers $n\leq x$ having property ($\star$)
such that $\chi(n)=\zeta$.  Then
$$
  \#\{n<x\mid \text{$n$ has property ($\star$)}\,,\; (n,m)=1\}
  \;=\;
  -
  \sum_{k=1}^{\ell-1}
  \zeta^{-k}
  \sumstar_{n\leq x}
  \chi^k(n)
  \,,
$$
where $\sum^\star$ means that the sum is taken only over those positive integers having property ($\star$).
\end{lemma}


\noindent\textbf{Proof.}
Summing the identity
$$
  \sum_{k=1}^\ell \zeta^{-k}\chi^k(n)=
  \begin{cases}
  \ell & \chi(n)=\zeta\\
  0 & \text{otherwise}
  \end{cases}
  \,.
$$
over all $n\leq x$ satisfying ($\star$)
and isolating the trivial character from the resulting expression
gives the desired conclusion.
\qed

\begin{lemma}\label{L:PV}
  Let $\chi$ be a non-principal Dirichlet character modulo $m\geq 2\cdot 10^4$, and
  let $p$ be a prime.
  For $x>0$, we have
  $$
  \vrule \
  \sum_{\substack{n<x\\(n,p)=1}}^{\phantom{n}}
  \chi(n)
  \ \vrule
  \leq
  2\sqrt{m}\log m
  \,.
  $$
\end{lemma}

\noindent\textbf{Proof.}
Given that $m\geq 2\cdot 10^4$,
the explicit version
of the P\'olya--Vinogradov inequality
given in~\cite{bachman.rachakonda}
implies that for
for any $y>0$, we have
\begin{equation}\label{E:PVproof1}
  \left|
  \sum_{n<y}
  \chi(n)
  \right|
  \leq
  m^{1/2}\log m
  \,.
\end{equation}
We write
\begin{equation}\label{E:PVproof2}
    \sum_{\substack{n<x\\(n,p)=1}}
  \chi(n)
  =
  \sum_{n<x}
  \chi(n)
  -
  \chi(p)\sum_{n<x/p}
  \chi(n)
  \,.
\end{equation}
Applying the triangle inequality to (\ref{E:PVproof2}) and invoking (\ref{E:PVproof1}) twice
gives the result.~\qed


\begin{lemma}\label{L:DC2}
Suppose $\chi$ is a Dirichlet character modulo $m$.
Suppose $q\geq 3$ is a positive integer, and let $A$ be a subset of $(\Z/q\Z)^*$.
Let $(\star)$ be any property of integers.
We have
$$
  \vrule \
  \sum_{a\in A}
  \;\;
  \sumstar_{\substack{n\leq x\\n\equiv a\pod{q}}}^{\phantom{n}}
  \chi(n)
  \ \vrule
  \leq
  \frac{\phi(q)}{2}
  \;
  \max_{\psimodq}
  \left|
  \sumstar_{n\leq x}
  \left(\psi\chi\right)(n)
  \right|
  \,,
$$
where $\sum^\star$ means that the sum is only taken over those positive integers $n$ having property $(\star)$.
\end{lemma}

\noindent\textbf{Proof.}
For notational convenience we denote $N:=\# A$.
We begin by summing the identity
$$
  \frac{1}{\phi(q)}
  \sum_{\psimodq}
  \overline{\psi}(a)\psi(n)\chi(n)
  =
  \begin{cases}
    \chi(n) & n\equiv a\pod{q} \\
    0 & \text{otherwise}
  \end{cases}
  \,,
$$
over all $n\leq x$ satisfying $(\star)$ and all $a\in A$, to obtain
\begin{eqnarray*}
      \sum_{a\in A}
  \;\;
    \sumstar_{\substack{n\leq x\\n\equiv a\pod{q}}}
    \chi(n)
    &=&
    \sum_{a\in A}
    \sumstar_{n\leq x}
    \frac{1}{\phi(q)}
    \sum_{\psimodq}
    \overline{\psi}(a)\psi(n)\chi(n)
    \\
    &=&
    \frac{1}{\phi(q)}
    \sum_{\psimodq}
    \left(
    \sum_{a\in A}
    \overline{\psi}(a)
    \right)
    \left(
    \sumstar_{n\leq x}
    \left(\psi\chi\right)(n)
    \right)
    \,.
\end{eqnarray*}

Observe that
$$
  \left|
    \sum_{a\in A}\overline{\psi}(a)
  \right|
  \leq
  N
  \,.
$$
Therefore if $N\leq \phi(q)/2$, we are done.
Hence we may assume that $\phi(q)/2+1\leq N\leq\phi(q)$.
In this case, we observe that
when $\psi$ is not the trivial character modulo $q$ we have
$$
  \sum_{a\in A}\overline{\psi}(a)
  =
  -
  \sum_{a\notin A}\overline{\psi}(a)
  \,,
$$
and the result follows upon observing that
\begin{eqnarray*}
  \frac{1}{\phi(q)}
  \sum_{\psimodq}
  \left|
  \sum_{a\in A}\overline{\psi}(a)
  \right|
  &\leq&
  \frac{(\phi(q)-N)(\phi(q)-1)+N}{\phi(q)}
  \\
  &\leq&
  \frac{\phi(q)}{2}
  \,.
  \;\;
  \text{\qed}
\end{eqnarray*}

\vspace{2ex}
\noindent\textbf{Proof of Proposition~\ref{P:special}.}
We may assume $f>\ell^2\geq 9$ as this is implied by either
inequality appearing in our hypothesis.
Now by Lemma~\ref{L:disc}, we may assume that
$f$ is a prime with $f\equiv 1\pmod{\ell}$.

First suppose that $q_1=2$ and $q_2=3$.
We will say that $n\in\Z^+$ has property $(\star)$ if
$(6,n)=1$ and $n\not\equiv 3f\pmod{4}$.
By condition~(1) of Theorem~\ref{T:conditions}, we must prove that there exists $r\in\Z^+$
satisfying condition $(\star)$ with $\chi(r)=\chi(3)^{-1}=:\zeta$ such that
$3r-1\leq f$.
By way of contradiction, suppose there are no positive integers $n<x$
satisfying condition $(\star)$ with $\chi(n)=\zeta$.
We will choose $x$ later, but for now, we assume $0<x<f$.

Applying Lemma~\ref{L:DC1}
we have:
\begin{equation}\label{E:q1is2}
  \#\{n<x\mid n \text{ has property $(\star)$} \}
  \;\leq\;
  (\ell-1)\max_{k=1,\dots,\ell-1}
  \left|
  \sumstar_{n<x}
  \chi^k(n)
  \right|
\end{equation}
First we estimate the quantity on the left-hand side of (\ref{E:q1is2}) from below.
Observe that:
\begin{eqnarray*}
    \#\{n<x\mid n \text{ has property $(\star)$} \}
    &=&
    \#\{n<x\;\mid\; n\,\equiv\, 3f+2,\, 3f+10\pmod{12}\}
    \\[1ex]
    &\geq&
    \frac{x}{6}-2
\end{eqnarray*}
Now we estimate the sum on the right-hand side of (\ref{E:q1is2}) from above.
By
Lemma~\ref{L:PV}
and
Lemma~\ref{L:DC2}, we have
\begin{eqnarray*}
  \left|
  \sumstar_{n<x}
  \chi^k(n)
  \right|
  &\leq&
  \max_{\psi\hspace{-1.25ex}\mod 4}
  \;
  \vrule \
    \sum_{\substack{n<x\\(3,n)=1}}^{\phantom{n}}
    (\psi\chi^k)(n)
 \ \vrule
  \\[1ex]
  &\leq&
  2(4f)^{1/2}\log 4f
  \,.
\end{eqnarray*}

Putting everything together, we have
$$
  \frac{x}{6}-2
  <
  4(\ell-1)f^{1/2}\log 4f
  \,,
$$
which implies
$$
   x<24(\ell-1)f^{1/2}\log 4f+12
   \,.
$$
Hence there exists an $r\in\Z^+$ with $\chi(r)=\zeta$ and
$$
  r\leq 24(\ell-1)f^{1/2}\log 4f+12
  \,,
$$
lest we arrive at a contradiction.
In light of this, to satisfy condition~(1) of Theorem~\ref{T:conditions},
which reads $3r-1\leq f$ in this case,
it is enough to assume
$$
  3(24(\ell-1)f^{1/2}\log 4f+12)-1\leq f
  \,,
$$
which is true by hypothesis.

Now we treat the second case of $q_1=3$ and $q_2=5$.
We only sketch the proof as it is very similar to the first.
This time, we will say that $n\in\Z^+$ has property $(\star)$ if
$(15,n)=1$ and $n\not\equiv f,2f\pmod{9}$;
we find that this holds exactly when $n$ belongs to one of
$16$ particular residue classes modulo $45$.
By condition~(1) of Theorem~\ref{T:conditions}, we must prove that there exists $r\in\Z^+$
satisfying condition $(\star)$ with $\chi(r)=\chi(5)^{-1}=:\zeta$ such that
$10r-2\leq f$.
By way of contradiction, suppose there are no positive integers $n<x$
satisfying condition $(\star)$ with $\chi(n)=\zeta$.

We find
\begin{eqnarray*}
    \#\{n<x\mid n \text{ has property $(\star)$} \}
    &>&
    \frac{16\,x}{45}-16
\end{eqnarray*}
and
\begin{eqnarray*}
  \left|
  \sumstar_{n<x}
  \chi^k(n)
  \right|
  &\leq&
  3\max_{\psi\hspace{-1.5ex}\mod 9}
  \;\;
   \vrule \
    \sum_{\substack{n<x\\(5,n)=1}}^{\phantom{n}}
    (\psi\chi^k)(n)
  \ \vrule
  \\[1ex]
  &\leq&
  6(9f)^{1/2}\log 9f
  \,.
\end{eqnarray*}
Combining the above, using the same argument as before,
we find
$$
  \frac{16}{45}\,x
  <
  18(\ell-1)f^{1/2}\log 9f+16
  \,.
$$
Proceeding as before, we arrive at the desired result.
\qed


\subsection{Upper bounds on $q_1$, $q_2$, and $r$}\label{S:DB.2}


Here we give bounds on the quantities $q_1$, $q_2$, and $r$ appearing in Theorem~\ref{T:conditions}.
First we quote the following result, which is proved elsewhere:
\begin{theorem}\label{T:small}
Let $\chi$ be a non-principal Dirichlet character modulo a prime $p\geq 10^{19}$
having odd order.
Suppose that $q_1<q_2$ are the two smallest prime non-residues of $\chi$.
Then we have:
\begin{enumerate}
  \item
  $q_1 < 3.9\, p^{1/4}\log p$
  \item
  $q_2<53\,p^{1/4}(\log p)^2$
  \item
  $q_1 q_2<24\,p^{1/2}(\log p)^2$
\end{enumerate}
\end{theorem}
The bound on $q_1$ above is due to Norton (see~\cite{norton:1971}), and the bounds
on $q_2$ and the product $q_1 q_2$ are due to the author (see~\cite{mcgown:small}).
In order to bound $r$, we will use the character sum estimate
(see Theorem~\ref{T:burgess}) given in the appendix; however, we remark that
Theorem~\ref{T:small} gives a stronger bound for $q_1$ and $q_2$ than one would achieve via Theorem~\ref{T:burgess}.

Now we state and prove a result which gives an upper bound on $r$.
Having dealt with the two special cases in \S\ref{S:DB.1},
we do not need to impose any additional
congruence conditions on $r$.
Larger values of $q_1$ lead to better constants, and so we provide two sets of constants.
\begin{proposition}\label{P:boundr}
  Let $\chi$ be a non-principal Dirichlet character modulo $f$ of order $\ell>2$,
  where $f$ is a prime with $f\geq 2\cdot 10^4$.
  Let $q_1<q_2$ be primes.
  Fix an $\ell$-th root of unity $\zeta$, and
  $k\in\Z$ with $k\geq 2$.
  There exists a computable positive constant $D(k)$ such that
  whenever $f$ is large enough so that
  $$
    \left(D(k)(\ell-1)\right)^k(\log f)^{\frac{1}{2}}
    \leq
    4f^{\frac{1}{4}}
    \,,
  $$
  there exists $r\in\Z^+$ such that $(r,q_1 q_2)=1$, $\chi(r)=\zeta$, and
  $$
    r\leq \left(D(k)\,(\ell-1)\right)^k\,f^{\frac{k+1}{4k}}(\log f)^{\frac{1}{2}}
    \,.
  $$

\begin{minipage}{2.5in}
\begin{table}[H]
\centering
\begin{tabular}{c c}
\begin{tabular}{| l | r |}
\hline
$k$ & $D_1(k)\;$ \\
\hline
$2$ & $89.1550$\\
$3$ & $43.1104$\\
$4$ & $31.9985$\\
$5$ & $26.9751$\\
$6$ & $24.1129$\\
$7$ & $22.2635$\\
$8$ & $20.9692$\\
\hline
\end{tabular}
\qquad
\begin{tabular}{| l | r |}
\hline
$k$ & $D_1(k)\;$ \\
\hline
$9$ & $20.0133$\\
$10$ & $19.2768$\\
$11$ & $18.6920$\\
$12$ & $18.2160$\\
$13$ & $17.8211$\\
$14$ & $17.4877$\\
$15$ & $17.2028$\\
\hline
\end{tabular}
\end{tabular}
\caption{Values of $D(k)$ when $2\leq k\leq 15$, with $q_1$ arbitrary\label{Tab:D1}}
\end{table}
\end{minipage}
\begin{minipage}{2.5in}
\begin{table}[H]
\centering
\begin{tabular}{c c}
\begin{tabular}{| l | r |}
\hline
$k$ & $D_2(k)\;$ \\
\hline
$2$ & $13.5958$\\
$3$ & $6.6415$\\
$4$ & $5.0420$\\
$5$ & $4.3220$\\
$6$ & $3.9103$\\
$7$ & $3.6430$\\
$8$ & $3.4550$\\
\hline
\end{tabular}
\qquad
\begin{tabular}{| l | r |}
\hline
$k$ & $D_2(k)\;$ \\
\hline
$9$ & $3.3154$\\
$10$ & $3.2075$\\
$11$ & $3.1215$\\
$12$ & $3.0513$\\
$13$ & $2.9929$\\
$14$ & $2.9434$\\
$15$ & $2.9011$\\
\hline
\end{tabular}
\end{tabular}
\caption{Values of $D(k)$ when $2\leq k\leq 15$, assuming $q_1>100$\label{Tab:D2}}
\end{table}
\end{minipage}
\end{proposition}

\noindent\textbf{Proof.}
Define the constant $C(k)$ as in Theorem~\ref{T:burgess},
and two more quantities which depend on $q_1$, $q_2$, $k$:
  $$
   K_1:=\left(1+q_1^{1/k-1}\right)\left(1+q_2^{1/k-1}\right)
    \,,
    \quad
    K_2:=\left(1-q_1^{-1}\right)\left(1-q_2^{-1}\right)
$$
Fix a constant $D(k)$ such that
  $$
  D(k)\geq
  \frac{K_1\left(1+C(k)^{-1}\right)}{K_2}\, C(k)
  \,.
  $$
We will show that the theorem holds for this choice of $D(k)$.
Set
$$
  x:=\left(D(k)(\ell-1)\right)^k f^{\frac{k+1}{4k}}(\log f)^{\frac{1}{2}}
  \,,
$$  
and
suppose there are no positive integers $n<x$
with $(n,q_1 q_2)=1$ and 
\mbox{$\chi(n)=\zeta$}.
We observe that $x\leq 4f^{\frac{1}{2}+\frac{1}{4k}}$ by hypothesis;
in particular, we find $x<4f^{5/8}<f$.

Applying Lemma~\ref{L:DC1}
we have:
\begin{equation}\label{E:q1arb}
  \#\{n<x\mid (n,q_1 q_2)=1 \}
  \;\leq\;
  (\ell-1)\max_{k=1,\dots,\ell-1}
  \left|
  \sum_{\substack{n<x\\(n,q_1 q_2)=1}}
  \chi^k(n)
  \right|
\end{equation}
We bound the left-hand side of (\ref{E:q1arb}) from below:
$$
  \#\{n<x\mid (n,q_1 q_2)=1 \}
  > (1-q_1^{-1})(1-q_2^{-1})x-2
$$
Now we wish to bound the character sum on right-hand side of (\ref{E:q1arb}) from above.
We fix an arbitrary $k\in\{1,\dots,\ell-1\}$, and
for notational convenience, we will write $\chi$ in place of $\chi^k$.
We have:
$$
  \sum_{\substack{n<x\\(n,q_1 q_2)=1}}
  \chi(n)
  =
  \sum_{n<x}\chi(n)
  -
  \chi(q_1)
  \sum_{n<x/q_1}\chi(n)
  -
  \chi(q_2)
  \sum_{n<x/q_2}\chi(n)
  +
   \chi(q_1 q_2)
  \sum_{n<x/q_1 q_2}\chi(n)
$$
Now we apply the triangle inequality to the above and invoke Theorem~\ref{T:burgess}
to bound each term.  This gives
$$
  \left|
  \sum_{\substack{n<x\\(n,q_1q_2)=1}}
  \chi^k(n)
  \right|
  \leq
  C(k)\left(1+q_1^{1/k-1}\right)\left(1+q_2^{1/k-1}\right)x^{1-1/k}f^{\frac{k+1}{4k^2}}(\log f)^{\frac{1}{2k}}
  \,.
$$  
  
Combining everything, we have
\begin{eqnarray*}
   K_2\,x
  &<&
  (\ell-1)
  K_1\,C(k)
  x^{1-\frac{1}{k}}f^{\frac{k+1}{4k^2}}(\log f)^{\frac{1}{2k}}
  +2
  \\
  &\leq&
(\ell-1)K_1\left(1+ C(k)^{-1}\right)C(k)x^{1-\frac{1}{k}}f^{\frac{k+1}{4k^2}}(\log f)^{\frac{1}{2k}}
  \,,
\end{eqnarray*}
which leads to
$$
  x
  <
  \left(D(k)(\ell-1)\right)^k f^{\frac{k+1}{4k}}(\log f)^{\frac{1}{2}}
  \,,
$$
a contradiction.\footnote{Computation of the table of constants is routine.  For the first set of constants, we use $q_1\geq 2$, $q_2\geq 3$, and for the second set we use $q_1\geq 101$, $q_2\geq 103$.}
\qed

\subsection{The general case}



%

Having paved the way, we are ready to prove the following result from which
Theorem~\ref{T:DB} follows immediately.

\begin{theorem}\label{T:DB.inequality}
  Fix an integer $2\leq k\leq 8$ and define the
  constant $E(k)$ as in Table~\ref{Table:E}.
  Let $K$ be a Galois number field of odd prime degree $\ell$ and conductor $f$.
  If
  $$
    E(k)(\ell-1)^k(\log f)^{\frac{7}{2}}\leq
    f^{\frac{1}{4}-\frac{1}{4k}}
    \,,
  $$
  then $K$ is not norm-Euclidean.

\begin{table}[H]
\centering
\begin{tabular}{c c}
\begin{tabular}{| l | r |}
\hline
$k$ & $E(k)\phantom{123}$ \\
\hline
$2$ & $3.4936\cdot 10^{3}$\\
$3$ & $5.5369\cdot 10^{3}$\\
$4$ & $1.2215\cdot 10^{4}$\\
$5$ & $2.8503\cdot 10^{4}$\\
$6$ & $6.7566\cdot 10^{4}$\\
$7$ & $1.6095\cdot 10^{5}$\\
$8$ & $3.8375\cdot 10^{5}$\\
\hline
\end{tabular}
\end{tabular}
\caption{Values of $E(k)$\label{Table:E}}
\end{table}
\end{theorem}

\vspace{1ex}
\noindent\textbf{Proof of Theorem~\ref{T:DB}.}
If $\ell=3$, then set $k=5$.  If $3<\ell<61$, then set $k=4$.  Otherwise set $k=3$.
Apply Theorem~\ref{T:DB.inequality}.\footnote{
Since any choice of $k$ will give a discriminant bound, we merely test numerically the values of $k\in[2,8]$ to see
which choice gives the least exponent in the bound.
It appears that after a certain point, $k=2$ will be the best choice.
}
\qed

\vspace{1ex}

We note in passing that we could derive a similar inequality to that given in
Theorem~\ref{T:DB.inequality} for all $k\geq 2$, but as these results will not
improve our ultimate discriminant bounds,
we have opted to use the simplifying assumption of $k\leq 8$.

\vspace{1ex}

\noindent\textbf{Proof of Theorem~\ref{T:DB.inequality}.}
Our ultimate choice of $E(k)$ will be such that $E(k)\geq 10^3$.  Using this, together with
$k\geq 2$, $\ell\geq 3$, our hypothesis leads to the inequality
\mbox{$4\cdot 10^3(\log f)^\frac{7}{2}\leq f^{\frac{1}{4}}$}
which easily implies
$f\geq 10^{40}$.
One also checks that \mbox{$f>\ell^2$}, from the hypothesis.
Using Lemma~\ref{L:disc}, we may assume $f$ is a prime with
\mbox{$f\equiv 1\pmod{\ell}$}.
We adopt the notation from the hypothesis of Theorem~\ref{T:conditions},
and set $\zeta=\chi(q_2)^{-1}$.

For now we will assume $q_1>100$.
Using Theorem~\ref{T:conditions}, we must show there exists
$r\in\Z^+$ such that $(r,q_1 q_2)=1$, $\chi(r)=\zeta$,
which also satisfies the inequality
$$
  2.1\, q_1 q_2 r\log q_1\leq f
  \,.
$$
Using Theorem~\ref{T:small}, we have
$$
  q_1 q_2<24\,f^{1/2}(\log f)^2
  \,,
$$
and
$
  q_1<3.9\,f^{1/4}\log f<f^{3/8}
$,
which implies
$$
  \log q_1<\frac{3}{8}\log f
  \,.
$$
Thus, we have
$$
  2.1\, q_1 q_2 \log q_1
  <
  18.9\, f^{1/2}(\log f)^3
  \,.
$$
Using Proposition~\ref{P:boundr} we obtain an integer $r$ with the desired properties such that
$$
  r\leq \left(D_2(k)\,(\ell-1)\right)^k\,f^{\frac{k+1}{4k}}(\log f)^{\frac{1}{2}}
  \,,
$$
provided
\begin{equation}\label{E:extracond}
    \left(D_2(k)(\ell-1)\right)^k(\log f)^{\frac{1}{2}}
    \leq
    4f^{\frac{1}{4}}
    \,.
\end{equation}

We define the constant
$$
  E(k):=18.9\, D_2(k)^k
  \,.
$$
Combining everything, and using the hypothesis,
we have the bound
$$
    2.1\, q_1 q_2 r\log q_1
    <
    E(k)(\ell-1)^k(\log f)^\frac{7}{2}f^{\frac{3k+1}{4k}}
    \leq
    f
    \,.
$$
It remains to verify (\ref{E:extracond}), but having defined $E(k)$,
we easily verify that this condition is automatic from our hypothesis
as one has:
\begin{eqnarray*}
  \left(D_2(k)(\ell-1)\right)^k(\log f)^{\frac{1}{2}}
  &\leq&
  E(k)(\ell-1)^k(\log f)^{\frac{1}{2}}
  \\
  &\leq&
   \frac{f^{\frac{1}{4}-\frac{1}{4k}}}{(\log f)^3}
   \\
   &<&
   f^{\frac{1}{4}}
\end{eqnarray*}
This completes the proof in the case that $q_1>100$.

Now we consider what happens when $q_1\leq 100$.
Having dealt with two special cases in \S\ref{S:DB.1},
the remaining cases fall under conditions~(2) through (5) of
Theorem~\ref{T:conditions}.
Namely, we must show there exists
$r\in\Z^+$ such that $(r,q_1 q_2)=1$, $\chi(r)=\zeta$,
which also satisfies an additional inequality.
We will prove the bound
\begin{equation}\label{E:q100want}
  932\, q_2 r<f
  \,,
\end{equation}
which will establish the result in all cases;
in particular, we observe that
$$
  (2.1)(97)(\log 97)<932
  \,.
$$  

We apply Lemma 7 and Theorem 4 of~\cite{mcgown:small} to find that $q_1\leq 100$ implies $q_2<711\,p^{1/4}\log p$.
Using Proposition~\ref{P:boundr} we obtain an integer $r$ with the desired properties such that
$$
  r\leq \left(D_1(k)\,(\ell-1)\right)^k\,f^{\frac{k+1}{4k}}(\log f)^{\frac{1}{2}}
  \,,
$$
provided
\begin{equation}\label{E:extracondprime}
    \left(D_1(k)(\ell-1)\right)^k(\log f)^{\frac{1}{2}}
    \leq
    4f^{\frac{1}{4}}
    \,.
\end{equation}
We obtain
$$
  932\, q_2 r<E'(k)(\ell-1)^k f^{\frac{1}{2}+\frac{1}{4k}}(\log f)^\frac{3}{2}
  \,,
$$
where
$$
  E'(k)=(932)(711)D_1(k)^k
  \,.
$$
To complete the proof, it suffices to show
\begin{equation}\label{E:finalsuffices}
E'(k)(\ell-1)^k f^{\frac{1}{2}+\frac{1}{4k}}(\log f)^\frac{3}{2}\leq f
\,,
\end{equation}
as (\ref{E:finalsuffices}) implies both (\ref{E:q100want}) and (\ref{E:extracondprime}).

But one checks that (\ref{E:finalsuffices})
follows from our hypothesis provided
\begin{equation}\label{E:q1smallcond}
  \frac{E'(k)}{E(k)}\leq f^{1/4}(\log f)^2
  \,.
\end{equation}
Finally, using the fact that $f\geq 10^{40}$ implies $f^{1/4}(\log f)^2\geq 10^{13}$,
an easy numerical computation shows that
(\ref{E:q1smallcond})
holds for $k=2,\dots,8$.
\qed

\section{An Algorithm and Some Computations}\label{S:comp}

In this section we give the algorithm to which we alluded in \S\ref{S:intro}.
In \S\ref{S:algorithm.idea} we give the main idea behind the algorithm,
in \S\ref{S:algorithm.character} we discuss character evaluations,
and in \S\ref{S:algorithm.statement} we give a full statement of the algorithm.
Finally, in \S\ref{S:algorithm.results} we give some results obtained from our computations,
including Theorems~\ref{T:comp} and~\ref{T:comp.cubic}.

\subsection{Idea behind the algorithm}\label{S:algorithm.idea}
Let us first state our aims in designing such an algorithm.
The input should be an odd prime $\ell$ and positive integers $A,B$.
If we let $\mathcal{F}_\ell(A,B)$ denote the collection of all Galois number fields $K$
of degree $\ell$ with conductor $f\in[A,B]$, then
the output should be a list $\mathcal{L}\subset[A,B]$
which contains the conductors of all norm-Euclidean
$K\in\mathcal{F}_\ell(A,B)$.

We do not require our list to consist of only norm-Euclidean fields, but 
the list should be manageable in the sense that we could eventually
hope to treat the remaining fields on a case-by-case basis.
Our goal is to sift through a very large amount of fields as quickly as possible.
We will use the first condition from Theorem~\ref{T:conditions} exclusively.
For the reader's convenience, we give this part of the theorem again:

\begin{theorem*}
  Let $K$ be a Galois number field of odd prime degree $\ell$ and conductor $f$ with $(f,\ell)=1$,
  and
  let $\chi$ be a primitive Dirichlet character modulo $f$ of order $\ell$.
  Denote by $q_1<q_ 2$ the two smallest rational primes with $\chi(q_1),\chi(q_2)\neq 1$.
  Suppose that there exists $r\in\Z^+$ with
  \begin{eqnarray*}
    &&
    \qquad
    (r,q_1 q_2)=1,
    \quad
    \chi(r)=\chi(q_2)^{-1},
    \\
    &&
    \qquad
    r q_2 k\not\equiv f \pmod{q_1^2}, \quad k=1,\dots,q_1-1\,,
    \\
    &&
    \qquad
    (q_1-1)(q_2 r-1)\leq f
    \,.
    \end{eqnarray*}
   Then $K$ is not norm-Euclidean.
\end{theorem*}

Although this is the most awkward condition (among those given in Theorem~\ref{T:conditions}) to apply in order to obtain theoretical bounds,
it is useful in computation as the congruence condition is satisfied more than half the time.
Indeed, if we assume the congruence class of $r$ inside $(\Z/ q_1^2\Z)^\star$ is chosen randomly, the chances
the condition is satisfied are $(q_1-1)/q_1$.
Therefore, when $q_1$ is large, it is very likely that any value we take for $r$ will automatically 
satisfy our congruences; on the other hand, when $q_1$ is small, the congruences may fail on occasion, but in this case we have lots of room to look
for $r$.  In addition, the conditions above only require computation within $\Z$ and character evaluations, and hence one can avoid
the additional considerations of precision that come along with computing logarithms.

In searching for the integer $r$ required to apply the above theorem, performing character evaluations is unavoidable.  The basic idea is to arrange things so that
character evaluations are almost the only computations needed, and that we carry out as few of them as possible.
To this end, our algorithm will only perform character
evaluations on primes.  This has the advantage that we won't have to sieve out a list of integers coprime to
$q_1 q_2$ for each $f$; instead, for each $f$
we evaluate a fixed character $\chi$ against a precomputed list of primes.

Based on the above discussion,
the basic strategy is as follows:
compute
$\chi(p)$ for primes $p<f$ until we find the smallest prime non-residues
$q_1$, $q_2$ and a prime $r$ with $\chi(r)=\chi(q_2)^{-1}$ satisfying our congruences.
If we are able to do this before we run out of primes, then we simply
check whether $(q_1-1)(q_2r-1)\leq f$.
Assuming any of the $\ell$-th roots of unity are equally likely to occur,
  and that our congruences are satisfied at least half the
  time,
  then an upper bound on the average
  number of character evaluations to find $q_1$, $q_2$, $r$ as just described is:
$$
  \ell\left(2+\frac{2}{\ell-1}\right)
$$

This gives a rough heuristic for how many character evaluations are necessary.
For example, when $\ell=3$, it should take almost 9 character evaluations on average
to prove that any given cubic field is not norm-Euclidean.\footnote{A quick test using the range $100\leq f\leq 300$ yields an average of $\approx 8.7$.}
However, it is important to keep in mind that on occasion it may take many more
character evaluations than the average.

\subsection{Character evaluations}\label{S:algorithm.character}
Before stating the algorithm formally, we detour for a brief discussion as to how we
carry out our character evaluations, as this
constitutes the largest portion of the computations.
Fix $\chi=\chi_f$, a primitive Dirichlet character modulo $f$ of order $\ell$ where $f$ and $\ell$ are both odd primes
and $f\equiv 1\pmod{\ell}$.

If we are performing multiple evaluations of a single character, and the modulus $f$ is small, then perhaps one of
the best strategies is to first build a lookup table.  Once this is completed, we can perform character evaluations in small constant time.
One straightforward way to do this is to first find a primitive root for $f$.  We won't go into algorithms for this here.

When $f$ is too large, building a lookup table is not a good option as it becomes infeasible to store
such a table in memory, and seems excessive given that it is very likely that we will only need to evaluate the
character a small number of times.
We describe an alternative approach based on the following observation:
If $\mathfrak{p}$ is a prime in $\Q(\zeta_\ell)$ with $\mathfrak{p}\divides f$,
then $\psi(n)=(n/\mathfrak{p})_\ell$ is a Dirichlet character modulo $f$ of order $\ell$,
where $(\,\cdot\,/\mathfrak{p})_\ell$ is the $\ell$-th power residue symbol (see \S4.1 of~\cite{lemmermeyer:book}).
In fact, the $\ell-1$ choices of $\mathfrak{p}$ lying over $f$ account for the $\ell-1$ Dirichlet characters modulo $f$ of order $\ell$;
hence we may assume
$\chi_f(n)=\left(n/\mathfrak{p}\right)_\ell$.

Moreover,
if we assume $\ell\leq 19$, then $\Q(\zeta_\ell)$ has class number one
(see~\cite{myron.montgomery})
and hence we may write
$\chi_f(n)=\left(n/\pi\right)_\ell$
for some prime $\pi\in\Z[\zeta_\ell]$ with $\pi\divides f$.\footnote
{
  The prime $\pi$ can be computed as $\gcd(\zeta_\ell-w, f)$, where $w$ is a solution
  to $\Phi_\ell(x)\equiv 0\pmod{f}$; here $\Phi_\ell(x)=x^{\ell-1}+\dots+x+1$ denotes the $\ell$-th cyclotomic polynomial.
}
In this case, we can use Eisenstein's reciprocity law for power residues to compute the above
symbol very rapidly, using computations in $\Z[\zeta_\ell]$,
in a manner completely analogous to the usual method of computing
Legendre symbols via the Jacobi symbol.
In this paper, we employ this  procedure in the cubic setting only;
see~\cite{cubicsymbol:2005} for details on how the computation of the cubic residue symbol can be carried out,
including the statement of the cubic reciprocity law.

\subsection{Statement of the algorithm}\label{S:algorithm.statement}

The input to our algorithm consists of positive integers $A,B$ and an odd prime $\ell$.
The output is a list $\mathcal{L}\subset[A,B]$ containing the conductors of all $K\in\mathcal{F}_\ell(A,B)$.
In the statement of Algorithm~\ref{A:1} below, a lowercase or uppercase latin letter will denote an element of $\Z$,
an uppercase script letter will denote a list of elements in $\Z$, and $\zeta$ will denote an $\ell$-th root of unity
(which can be stored as an integer in the interval $[0,\ell)$).  We will denote by $\chi_f$ a primitive Dirichlet character modulo $f$ of order $\ell$;
it does not matter which one we take as long as we use just one character for each $f$.


\begin{algorithm}\label{Alg:1}
\caption{Output a list of all possible conductors $f\in[A,B]$}
\label{A:1}
\begin{algorithmic}[1]
\STATE Generate a list $\mathcal{P}$ of all primes $p\leq\max\{1000, \sqrt{B}\}$ using the Sieve of Eratosthenes.
\STATE Generate a list $\mathcal{F}$ all primes $f\in[A,B]$ such that $f\equiv 1\pmod{\ell}$.
\FOR{$f\in\mathcal{F}$}
	\STATE Initiate scheme to evaluate $\chi_f$ (see~\S\ref{S:algorithm.character}).
	\STATE $q_1\gets 0$; $q_2\gets 0$; $r\gets0$
	\FOR{$p\in\mathcal{P}$}
		\IF{$p\geq f$}
			\STATE \textbf{break}
		\ENDIF
		\IF{\textbf{(}$\chi_f(p)\neq 1$\textbf{)}}
			\IF{$q_1=0$}
				\STATE $q_1\gets p$
			\ELSIF{$q_2=0$}
				\STATE $q_2\gets p$
				\STATE $\zeta\gets \chi_f(p)^{-1}$
				\STATE $\mathcal{A}\gets\{f q_2^{-1}k^{-1}\mod q_1^2\mid k=1,\dots,q_1-1\}$
			\ELSIF {$\chi_f(p)=\zeta$ AND $p\text{ mod }q_1^2\notin\mathcal{A}$}
			  \STATE $r\gets p$
			  \STATE \textbf{break}
			\ENDIF
			
		\ENDIF
	\ENDFOR
	\IF{$r= 0$ OR $(q_1-1)(q_2 r-1)>f$}
	  \PRINT $f$
	\ENDIF
\ENDFOR
\IF{$\ell^2\in[A,B]$}
	\PRINT $\ell^2$
\ENDIF

\end{algorithmic}
\end{algorithm}
Verifying the correctness of Algorithm~\ref{A:1} is straightforward.
For a given $f$, our algorithm either finds $q_1$, $q_2$, and $r$ satisfying the conditions in the theorem
or it doesn't; if it doesn't, then that value of $f$ is outputted.
However, we do give a number of comments regarding the algorithm which we
deem to be relevant:
  \begin{enumerate}
  \item
    In line 1, the reason for the number $1000$ is that if $B$ is especially small, we don't want to run out of primes.
    Of course, the number $1000$ is arbitrary -- any relatively manageable number will do.
   \item
   If we do run out of primes, the value of $r$ will remain at zero when the loop over $\mathcal{P}$ finishes.
   This will cause the relevant value of $f$ to be output, and so we need not worry about missing an $f$ due to lack of primes,
   or due to the non-existence of the value $r$ for that matter.
   \item
   In calculating the list $\mathcal{F}$ in line 2, one should sieve using the primes in $\mathcal{P}$ -- this is why we stored primes up to $\sqrt{B}$.
   \item
   Notice that the command ``Initiate scheme to evaluate $\chi_f$'' on line 4 is only run once for each $f$.
   Whether we are building a lookup table or finding a prime $\pi$ over $f$ (see~\S\ref{S:algorithm.character}),
   this step is carried out just once and results in fast character evaluations during the inner loop over $\mathcal{P}$.
   \item
   Although $\chi_f(p)$ appears on lines 10, 15, and 17, we of course only compute $\chi_f(p)$ once.
   \item
   The code on lines 15 and 16 to store values in $\zeta$ and $\mathcal{A}$ only gets executed at most once for each $f$.
   \item
   The modular arithmetic that takes place on lines 16 and 17 is modulo $q_1^2$, and typically $q_1$
   is very small.\footnote{Using rough heuristics as in \S\ref{S:algorithm.character}, we find that in the cubic case $q_1\in\{2,3,5,7\}$ roughly $98.8\%$ of the time,
   and as $\ell$ gets larger, this probability increases.}
  \end{enumerate}

\subsection{Results of the computations}\label{S:algorithm.results}

We have implemented the algorithm in the mathematics software SAGE\footnote{http://www.sagemath.org}
using a lookup table for character evaluations.
Using this, the results given in Theorem~\ref{T:comp}
took only 18.3 minutes of CPU time to complete
on a MacBook Pro with a 2.26 GHz Intel Core 2 Duo processor and 4 GB of RAM, running Mac OS 10.6.

For the cubic case, we have implemented an efficient version of our algorithm in C,
performing character evaluations using the equality $\chi_f(n)=(n/\pi)_3$,
as described in \S\ref{S:algorithm.character}.
We use NTL with GMP for large integer arithmetic, and 
we use the algorithms given in~\cite{cubicsymbol:2005} to compute the cubic residue symbol
and the greatest common divisor in $\Z[\zeta_3]$.
Running this code on all conductors $f\leq 10^{10}$ produced the same list of conductors as the $\ell=3$ entry in
Table~\ref{Table:candidate}.
This took
91.3 hours of CPU time on an iMac with a
3.06 GHz Intel Core 2 Duo processor and 4 GB of RAM, running Mac OS 10.6.
Thus we obtain:
\begin{theorem}\label{T:comp2}
There are no norm-Euclidean Galois cubic fields with discriminant
$1597^2<|\Delta|<10^{20}$.
\end{theorem}

Now we have:

\vspace{1ex}
\noindent\textbf{Proof of Theorem~\ref{T:comp.cubic}.}
Combine Theorems~\ref{T:godwin.smith}, \ref{T:DB}, and \ref{T:comp2}.
\qed

\vspace{1ex}

Not only does Theorem~\ref{T:comp2} extend the computations given in Theorem~\ref{T:comp},
but it provides a consistency check for the implementation of our character evaluations in both cases as the two implementations
are in two different programming languages using two completely different strategies for character evaluation.
We give the values of $q_1$, $q_2$, $r$
for the last 10 fields in our computation:

\begin{verbatim}
     f=9999999673, q1=5, q2=7, r=17
     f=9999999679, q1=2, q2=3, r=19
     f=9999999703, q1=2, q2=3, r=11
     f=9999999727, q1=7, q2=11, r=19
     f=9999999769, q1=3, q2=5, r=37
     f=9999999781, q1=2, q2=5, r=7
     f=9999999787, q1=3, q2=5, r=29
     f=9999999817, q1=2, q2=3, r=13
     f=9999999943, q1=5, q2=7, r=19
     f=9999999967, q1=5, q2=7, r=11
\end{verbatim}


\section*{Acknowledgments.}
This paper is based upon a portion of the author's Ph.D. dissertation under the supervision of Professor Harold Stark.  The author would like to thank
Professor Stark for his invaluable guidance and support at all stages of this work.

The author would also like to thank Professor Tom Schmidt for his comments on a preliminary version of this manuscript.



\bibliographystyle{amsplain}


\bibliography{myrefs}  

\providecommand{\bysame}{\leavevmode\hbox to3em{\hrulefill}\thinspace}
\providecommand{\MR}{\relax\ifhmode\unskip\space\fi MR }
\providecommand{\MRhref}[2]{%
  \href{http://www.ams.org/mathscinet-getitem?mr=#1}{#2}
}
\providecommand{\href}[2]{#2}
\begin{thebibliography}{10}

\bibitem{bachman.rachakonda}
G.~Bachman and L.~Rachakonda, \emph{On a problem of {D}obrowolski and
  {W}illiams and the {P}\'olya-{V}inogradov inequality}, Ramanujan J.
  \textbf{5} (2001), no.~1, 65--71. \MR{MR1829809 (2002c:11098)}

\bibitem{booker}
A.~Booker, \emph{Quadratic class numbers and character sums}, Math. Comp.
  \textbf{75} (2006), no.~255, 1481--1492 (electronic). \MR{MR2219039
  (2008a:11140)}

\bibitem{burgess:1962}
D.~A. Burgess, \emph{On character sums and primitive roots}, Proc. London Math.
  Soc. (3) \textbf{12} (1962), 179--192. \MR{MR0132732 (24 \#A2569)}

\bibitem{cubicsymbol:2005}
I.~Damg{\aa}rd and G.~Frandsen, \emph{Efficient algorithms for the gcd and
  cubic residuosity in the ring of {E}isenstein integers}, J. Symbolic Comput.
  \textbf{39} (2005), no.~6, 643--652. \MR{MR2168612 (2006g:11246)}

\bibitem{erdos.ko}
P.~Erd\"os and C.~Ko, \emph{Note on the {E}uclidean algorithm}, J. London Math.
  Soc. \textbf{13} (1938), 3--8.

\bibitem{friedlander:1984}
J.~Friedlander, \emph{Primes in arithmetic progressions and related topics},
  Analytic number theory and {D}iophantine problems ({S}tillwater, {OK}, 1984),
  Progr. Math., vol.~70, Birkh\"auser Boston, Boston, MA, 1987, pp.~125--134.
  \MR{MR1018373 (90h:11086)}

\bibitem{garbanati:cft}
D.~Garbanati, \emph{Class field theory summarized}, Rocky Mountain J. Math.
  \textbf{11} (1981), no.~2, 195--225. \MR{MR619671 (82g:12010)}

\bibitem{godwin:1965}
H.~J. Godwin, \emph{On the inhomogeneous minima of totally real cubic
  norm-forms}, J. London Math. Soc. \textbf{40} (1965), 623--627. \MR{MR0184927
  (32 \#2398)}

\bibitem{godwin.smith:1993}
H.~J. Godwin and J.~R. Smith, \emph{On the {E}uclidean nature of four cyclic
  cubic fields}, Math. Comp. \textbf{60} (1993), no.~201, 421--423.
  \MR{MR1149291 (93d:11114)}

\bibitem{hardy.wright}
G.~H. Hardy and E.~M. Wright, \emph{An introduction to the theory of numbers},
  fifth ed., The Clarendon Press Oxford University Press, New York, 1979.
  \MR{MR568909 (81i:10002)}

\bibitem{heilbronn:cyclic}
H.~Heilbronn, \emph{On {E}uclid's algorithm in cyclic fields}, Canadian J.
  Math. \textbf{3} (1951), 257--268.

\bibitem{ishida:book}
M.~Ishida, \emph{The genus fields of algebraic number fields}, Lecture Notes in
  Mathematics, Vol. 555, Springer-Verlag, Berlin, 1976. \MR{MR0435028 (55
  \#7990)}

\bibitem{iwaniec.kowalski}
H.~Iwaniec and E.~Kowalski, \emph{Analytic number theory}, American
  Mathematical Society Colloquium Publications, vol.~53, American Mathematical
  Society, Providence, RI, 2004. \MR{MR2061214 (2005h:11005)}

\bibitem{lemmermeyer:euclidean}
F.~Lemmermeyer, \emph{The {E}uclidean algorithm in algebraic number fields},
  Exposition. Math. \textbf{13} (1995), no.~5, 385--416. \MR{MR1362867
  (96i:11115)}

\bibitem{lemmermeyer:book}
\bysame, \emph{Reciprocity laws}, Springer Monographs in Mathematics,
  Springer-Verlag, Berlin, 2000, From Euler to Eisenstein. \MR{MR1761696
  (2001i:11009)}

\bibitem{myron.montgomery}
J.~M. Masley and H.~L. Montgomery, \emph{Cyclotomic fields with unique
  factorization}, J. Reine Angew. Math. \textbf{286/287} (1976), 248--256.
  \MR{MR0429824 (55 \#2834)}

\bibitem{mcgown:grh}
K.~McGown, \emph{Norm-{E}uclidean {G}alois fields and the {G}eneralized
  {R}iemann {H}ypothesis},  (in preparation).

\bibitem{mcgown:consecutive}
\bysame, \emph{On the constant in {B}urgess' bound for the number of
  consecutive residues or non-residues},  (submitted).

\bibitem{mcgown:small}
\bysame, \emph{On the second smallest prime non-residue},  (submitted).

\bibitem{norton:1971}
K.~Norton, \emph{Numbers with small prime factors, and the least {$k$}th power
  non-residue}, Memoirs of the American Mathematical Society, No. 106, American
  Mathematical Society, Providence, R.I., 1971. \MR{MR0286739 (44 \#3948)}

\bibitem{rosser.schoenfeld}
J.~B. Rosser and L.~Schoenfeld, \emph{Approximate formulas for some functions
  of prime numbers}, Illinois J. Math. \textbf{6} (1962), 64--94. \MR{MR0137689
  (25 \#1139)}

\bibitem{smith:1969}
J.~R. Smith, \emph{On {E}uclid's algorithm in some cyclic cubic fields}, J.
  London Math. Soc. \textbf{44} (1969), 577--582. \MR{MR0240075 (39 \#1429)}

\end{thebibliography}


\appendix
\section{An Explicit Version of Burgess' Character Sum Estimate}\label{S:burgess}

In this appendix, we prove an explicit version of a character sum estimate of Burgess
(see~\cite{burgess:1962}), following a method due to Iwaniec (see~\cite{iwaniec.kowalski} and~\cite{friedlander:1984}).
Booker proves a similar result when $\chi$ is quadratic (see~\cite{booker}).
\begin{theorem}\label{T:burgess}
Suppose $\chi$ is a non-principal Dirichlet character
modulo a prime $p\geq 2\cdot 10^4$.
Let $N,H\in\Z$ with $H\geq 1$.
Fix a positive integer $r\geq 2$.
Then there exists a computable constant $C(r)$ such that
whenever $H\leq 4p^{\frac{1}{2}+\frac{1}{4r}}$ we have
$$
  \left|
  \sum_{n\in(N,N+H]}
  \chi(n)
  \right|
  <
  C(r)\,
  H^{1-\frac{1}{r}}p^{\frac{r+1}{4r^2}}(\log p)^{\frac{1}{2r}}
  \,.
$$

\begin{table}[H]
\centering
\begin{tabular}{c c}
\begin{tabular}{| l | r |}
\hline
$r$ & $C(r)\phantom{1}$\\
\hline
$2$ & $10.0366$\\
$3$ & $4.9539$\\
$4$ & $3.6493$\\
$5$ & $3.0356$\\
$6$ & $2.6765$\\
$7$ & $2.4400$\\
$8$ & $2.2721$\\
\hline
\end{tabular}
\qquad
\begin{tabular}{| l | r |}
\hline
$r$ & $C(r)$\phantom{1}\\
\hline
$9$   & $2.1467$\\
$10$ & $2.0492$\\
$11$ & $1.9712$\\
$12$ & $1.9073$\\
$13$ & $1.8540$\\
$14$ & $1.8088$\\
$15$ & $1.7700$\\
\hline
\end{tabular}
\end{tabular}
\caption{Values for the constant $C(r)$ when $2\leq r\leq 15$:}
\end{table}
\end{theorem}
We note in passing that the assumption $H\leq 4p^{\frac{1}{2}+\frac{1}{4r}}$ is of a technical nature.
However, it seems that to drop it, at least in the current proof,
one may have to accept the slightly worse exponent of $1/r$ on the $\log p$ term.
In any case, this condition is essentially automatic for our application in \S\ref{S:DB}.

Throughout this section, $\chi$ will denote a Dirichlet character modulo an odd prime $p$
and $N,H$ will be integers with $0\leq N<p$ and $1\leq H<p$.
The latter assumption is justified as 
reducing $N$ and $H$ modulo $p$ leaves the sum in the above theorem unchanged.
The letter $r$ will denote a positive integer parameter with $r\geq 2$.
We begin with some definitions.

\begin{definition}
  $$
    S_\chi(H):=\sum_{n\in(N,N+H]}
    \chi(n)
   $$ 
\end{definition}
\begin{definition}
$$
  E(H):=H^{1-\frac{1}{r}}p^{\frac{r+1}{4r^2}}(\log p)^{\frac{1}{2r}}
$$  
\end{definition}
We seek a bound of the form
$S_\chi(H)<C\,E(H)$.
(An explicit way of choosing $C$ will appear in the statement of Theorem~\ref{T:charsum}.)
It is plain that $S_\chi(H)$ also depends upon $N$ and that $E(H)$ also depends upon $p$ and $r$,
but we have chosen to avoid excess decoration of our notations.  
\begin{definition}
  Fix $A\in\Z$ with $1<A<p$.
  For $x\in\F_p$, we define $\nu_A(x)$ 
  to be the number of ways we can write 
  $$
    x\equiv \overline{a}n\pmod{p}
    \,,
  $$  
  where $a\in[1,A]$ is a prime and $n\in(N,N+H]$ is an integer.
\end{definition}
In the above definition and in the rest of this section $\overline{a}$ will
denote a multiplicative inverse of $a$ modulo $p$.
We note that $\nu_A(x)$ also depends upon $N,H,p$.
Before launching the main part of the proof, we will require a series of lemmas.
%
%
\begin{lemma}\label{L:induct}
Suppose $|S_\chi(H_0)|\leq C\,E(H_0)$ for all $H_0<H$.
Fix $H_0=AB<H$.  Then
$$
  \left|S_\chi(H)\right|
  \leq
  \frac{1}{\pi(A)B}
  \sum_{x\in\F_p}
  \nu_A(x)
  \left|
  \sum_{1\leq b\leq B}
  \chi(x+b)
  \right|
  +
  2C\,E(H_0)
  \,.
$$
\end{lemma}

\noindent\textbf{Proof.}
Applying a shift $n\mapsto n+h$ with $1\leq h\leq H_0$ gives
\begin{equation}\label{E:shift}
\nonumber
  S_\chi(H)=\sum_{n\in(N,N+H]}
  \chi(n+h)
  +2C\theta E(H_0)
  \,.
\end{equation}
(The letter $\theta$ will denote a complex number with $|\theta|\leq 1$, possibly different each time it appears.)
We set $h=ab$ in the above, and average over all primes $a\in[1,A]$ and all integers $b\in[1,B]$.
This gives
$$
  S_\chi(H)
  =
  \frac{1}{\pi(A)B}
  \sum'_{a,b}
  \sum_{n\in(N,N+H]}
  \chi(n+ab)
  +
  2C\theta E(H_0)
  \,,
$$
where $\sum'$ here indicates that we are summing over all primes $a\in[1,A]$ and all integers $b\in[1,B]$.
Rearranging the sum in the above expression yields
$$
    \sum'_{a,b}
    \sum_{n\in(N,N+H]}
    \chi(n+ab)
    =
    \sum_{\substack{1\leq a\leq A\\\text{$a$ prime}}}
    \;
    \sum_{n\in(N,N+H]}
    \chi(a)
    \sum_{1\leq b\leq B}
    \chi(\overline{a}n+b)
    \,,
$$
and hence
$$
\left|
    \sum'_{a,b}
    \sum_{n\in(N,N+H]}
    \chi(n+ab)
\right|
\leq
\sum_{x\in\F_p}
\nu_A(x)
  \left|
  \sum_{1\leq b\leq B}
  \chi(x+b)
  \right|
  \,.
$$
The result follows.
\qed

%

\begin{lemma}\label{L:line}
  Suppose $a_1\neq a_2$ are prime and $b\in\Z$.
  Then the number of integer solutions $(x,y)\in\Z^2$ to the
  equation $a_1x - a_2 y =b$
  with $x,y\in(N,N+H]$ is at most
  $$
    \frac{H}{\max\{a_1,a_2\}}+1
    \,.
  $$  
\end{lemma}

\noindent\textbf{Proof.}
Let $Q$ denote the number of solutions to
  $a_1x - a_2 y =b$
with $x,y\in(N,N+H]$.
We will show
$Q\leq H/a_2 +1$.
It will immediately follow from the same argument that
$Q\leq H/a_1+1$
as well; indeed, just multiply both sides of the equation by $-1$
and interchange the roles of $x$ and $y$.
Suppose we have two solutions $(x,y)$ and $(x',y')$.
It follows that
$a_1(x-x')=a_2(y-y')$,
and since $a_1\neq a_2$ are prime, we see that $a_2$ divides $x-x'$ which implies $|x-x'|\geq a_2$.
The result follows.
\qed



\begin{lemma}\label{L:secondsum}
  Fix $A\in\Z$ with $1<A<p$.  If $2AH\leq p$, then
  \begin{eqnarray*}
    &&
    \sum_{x\in\F_p}
    \nu_A(x)^2
    <
    \pi(A)H
    \left(1+\frac{2}{\pi(A)}
    \hspace{-1ex}
     \sum_{\substack{a\leq A\\\text{$a$ prime}}}
    \frac{\pi(a)-1}{a}
    +
    \frac{2}{\pi(A)H}
     \hspace{-1ex}
    \sum_{\substack{a\leq A\\\text{$a$ prime}}}
    (\pi(a)-1)
    \right)
    \,.
  \end{eqnarray*}
\end{lemma}

\noindent\textbf{Proof.}
Define $S$ to be the set of all 
quadruples $(a_1,a_2,n_1,n_2)$
with
$$
  a_1 n_2\equiv a_2 n_1\pmod{p}
$$
where $a_1,a_2\in[1,A]$ are prime and $n_1,n_2\in(N,N+H]$ are integers.
We observe that
$\#S=\sum_{x\in\F_p}\nu_A(x)^2$.
Suppose $(a_1,a_2,n_1,n_2)\in S$ with $a_1=a_2$.
Then we have $n_1\equiv n_2\pmod {p}$ and hence $n_1=n_2$ since $n_1,n_2\in(N,N+H]$ and $H\leq p$.
Thus there are exactly $\pi(A)H$ solutions of this form.

Now we treat the remaining cases.
Let $(a_1,a_2,n_1,n_2)\in S$ with $a_1\neq a_2$.  Then
$a_1n_2-a_2n_1=kp$ for some $k$.
Writing $n_1=N+h_1$ and $n_2=N+h_2$ with $0<h_1,h_2\leq H$,
we have
\begin{eqnarray*}
k
&=&
\frac{a_1(N+h_2)-a_2(N+h_1)}{p}
\\[1ex]
&=&
\frac{a_1-a_2}{p}\,N
+
\frac{a_1 h_2-a_2 h_1}{p}
\\[1ex]
&=&
\frac{a_1-a_2}{p}\left(N+\frac{H}{2}\right)
+
\frac{a_1(h_2-H/2)-a_2(h_1-H/2)}{p}
\,,
\end{eqnarray*}
which gives
\begin{eqnarray*}
  \left|
  k-\left(\frac{a_1-a_2}{p}\right)\left(N+\frac{H}{2}\right)
  \right|
  <
  \frac{(a_1+a_2)H}{2p}
  \leq
  \frac{AH}{p}
  \leq
  \frac{1}{2}
  \,.
\end{eqnarray*}
This implies that $a_1$ and $a_2$ determine $k$.
Now Lemma~\ref{L:line} tells us that there are at most
$$\frac{H}{\max\{a_1,a_2\}}+1
$$
choices of $(n_1,n_2)$
for each fixed $(a_1,a_2)$.  Thus the number of elements in $S$ with $a_1\neq a_2$ is bounded above by
\begin{eqnarray*}
 2
 \sum_{\substack{a_2\leq A\\\text{$a_2$ prime}}}
 \;
 \sum_{\substack{a_1< a_2\\\text{$a_1$ prime}}}
  \left(
  \frac{H}{a_2}+1
  \right)
  &<&
  2H
   \sum_{\substack{a\leq A\\\text{$a$ prime}}}
  \frac{\pi(a)-1}{a}
  +
  2
   \sum_{\substack{a\leq A\\\text{$a$ prime}}}
  (\pi(a)-1)
  \,.
\end{eqnarray*}
This gives the result.
\qed

\vspace{1ex}
The next estimate is very weak, but has the advantage that it holds for all~$X$.
\begin{lemma}\label{L:pisum}
For $X\in\Z^+$ we have
$$
  \frac{1}{\pi(X)}
  \sum_{\substack{a\leq X\\a\text{ prime}}}
  \frac{\pi(a)-1}{a}
  <
  \frac{1}{3}
  \,.
$$
\end{lemma}

\noindent\textbf{Proof.}
The result holds for $X\leq 100$ by direct computation.
Using the Sieve of Eratosthenes, one easily shows that
$$
  \frac{\pi(n)-1}{n}\leq \frac{1}{3}
$$
for all $n\geq 100$.
The result follows.
\qed

\vspace{1ex}
Now we are ready to state and prove what is essentially the main result of this appendix,
from which Theorem~\ref{T:burgess} follows.

\begin{theorem}\label{T:charsum}
Suppose $\chi$ is a non-principal Dirichlet character
modulo an odd prime $p$.
Fix a positive integer $r\geq 2$.
Suppose $d>4$, $C\geq 1$, $p_0\geq 2$ are real constants satisfying
\begin{equation}\label{E:charsum.H1}
  C^rp_0^{\frac{1}{4}-\frac{1}{4r}}(\log p_0)^{\frac{1}{2}}
  \geq
  4d(d+1)r
\end{equation}
and
\begin{equation}\label{E:charsum.H2}
  C
  \geq
  \frac{
  \left(
  (d+1)(2r-1)(4r-1)
  \right)^{\frac{1}{2r}}
  }
  {
    \left(1-  \frac{2}{d^{1-\frac{1}{r}}}\right)
  }
  \,.
\end{equation}
If
$$
  H\leq \sqrt{rd}\,p^{\frac{1}{2}+\frac{1}{4r}}
  \,,
$$
then for $p\geq p_0$ we have
  $$
      |S_\chi(H)|\leq C\, E(H)
    \,.
  $$  
\end{theorem}

\noindent\textbf{Proof.}
We may assume
$$
  H\geq C^r p^{\frac{1}{4}+\frac{1}{4r}}(\log p)^{\frac{1}{2}}
  \,,
$$ 
or else the result follows from
the trivial bound $|S_\chi(H)|\leq H$.
We will prove the result by induction on $H$.
We assume that $|S_\chi(H_0)|\leq CE(H_0)$ for all $H_0<H$.
We choose an integer $H_0$ with
$$
  \frac{H}{d+1}<H_0\leq \frac{H}{d}
  \,,
$$
for which we can write
$H_0=AB$ with $A,B\in\Z^+$, where
$$
  B=\lfloor 4r p^{\frac{1}{2r}}\rfloor
  \,.
$$  
Accomplishing this is possible provided
$$
  H\geq 4d(d+1)rp^{\frac{1}{2r}}
  \,;
$$ 
given our a priori lower bound on $H$,
this condition follows from~(\ref{E:charsum.H1}).

Before proceeding further, we give upper and lower bounds on $A$.
Observe that
$$
  A\leq
  \frac{H}{dB}
  \leq
  \frac{\sqrt{rd}\,p^{\frac{1}{2}+\frac{1}{4r}}}{3drp^{\frac{1}{2r}}}
  =
  \frac{1}{3\sqrt{rd}}p^{\frac{1}{2}-\frac{1}{4r}}
  \,.
$$
We also have
$$
  A
  >
  \frac{H}{(d+1)B}
  \geq
  \frac{C^rp^{\frac{1}{4}+\frac{1}{4r}}(\log p)^{\frac{1}{2}}}{(d+1)4rp^{\frac{1}{2r}}}
  =
  \frac{C^rp^{\frac{1}{4}-\frac{1}{4r}}(\log p)^{\frac{1}{2}}}{4(d+1)r}
  \,.
$$
In particular, using
(\ref{E:charsum.H1}),
we see that
$A>d>4$.

Applying Lemma~\ref{L:induct} and our inductive hypothesis, we have
\begin{eqnarray}
\nonumber
  \left|S_\chi(H)\right|
  &\leq&
  \frac{1}{\pi(A)B}
  \sum_{x\in\F_p}
  \nu_A(x)
  \left|
  \sum_{1\leq b\leq B}
  \chi(x+b)
  \right|
  +
  2C\,E(H_0)
  \\
  \label{E:charsum.finally}
  &\leq&
    \frac{1}{\pi(A)B}
  \sum_{x\in\F_p}
  \nu_A(x)
  \left|
  \sum_{1\leq b\leq B}
  \chi(x+b)
  \right|
+
  \frac{2C}{d^{1-\frac{1}{r}}}E(H)
\,.
\end{eqnarray}  
In order to bound the sum above, we apply
H\"older's inequality to the functions
$\nu_A(x)^{1-\frac{1}{r}}$, $\nu_A(x)^\frac{1}{r}$, and $\left|\sum_{1\leq b\leq B}\chi(x+b)\right|$
using the H\"older exponents \mbox{$(1-1/r)^{-1}$}, $2r$, and $2r$ respectively; this yields:
\begin{eqnarray*}
  &&
  \sum_{x\in\F_p}
  \nu_A(x)
  \left|
  \sum_{1\leq b\leq B}
  \chi(x+b)
  \right|
  \\
  &&\qquad
  \leq
  \left(
  \sum_{x\in\F_p}
  \nu_A(x)
  \right)^{1-\frac{1}{r}}
    \left(
  \sum_{x\in\F_p}
  \nu_A(x)^2
  \right)^{\frac{1}{2r}}
  \left(
  \sum_{x\in\F_p}
  \left|
  \sum_{1\leq b\leq B}
  \chi(x+b)
  \right|^{2r}
  \right)^{\frac{1}{2r}}
\end{eqnarray*}

We bound each of the three sums above in turn.  Clearly, one has
$$
  \sum_{x\in \F_p}\nu_A(x)
  =\pi(A)H
  \,.
$$
We will shortly apply Lemma~\ref{L:secondsum} to show that
\begin{equation}\label{E:secondsumclaim}
  \sum_{x\in\F_p}
  \nu_A(x)^2
  \leq
  2\pi(A)H
  \,,
\end{equation}
but first we need to make a few estimates which involve the relevant quantities.

Our upper bound on $A$
allows us to verify that $2AH<p$, which makes Lemma~\ref{L:secondsum} applicable.
Lemma~\ref{L:pisum} gives
$$
  \frac{2}{\pi(A)}
  \sum_{\substack{a\leq A\\a\text{ prime}}}
  \frac{\pi(a)-1}{a}
  <
  \frac{2}{3}
  \,.
$$
Using (3.6) of~\cite{rosser.schoenfeld}, we have $\pi(A)\leq 1.26A/\log A$ for $A>1$ and therefore
\begin{eqnarray*}
\frac{\pi(A)}{H}
\leq
\frac{1.26 A}{H\log A}
\leq
\frac{1.26}{dB\log A}
\leq
\frac{1.26}{d(4r-1)\log A}
\leq
\frac{1.26}{4(4\cdot 2-1)\log 4}
<
0.1
\,.
\end{eqnarray*}
Now we see that
$$
  \frac{2}{\pi(A)H}
  \sum_{\substack{a\leq A\\a\text{ prime}}}
  (\pi(a)-1)
  \leq
  \frac{2\pi(A)}{H}
  <0.2
  \,.
$$
Putting all this together, we have successfully verified (\ref{E:secondsumclaim})
by invoking Lemma~\ref{L:secondsum}.

To bound the third sum, we apply Lemma~2.2 of~\cite{mcgown:consecutive}; this gives
$$
  \sum_{x\in\F_p}
  \left|
  \sum_{1\leq b\leq B}
  \chi(x+b)
  \right|^{2r}
  \leq
  B^{2r}p^{1/2}\left[\frac{1}{4}\left(\frac{4r}{B}\right)^rp^{1/2}+(2r-1)\right]
  \,.
$$
Notice that $B+1>4rp^{\frac{1}{2r}}$, and, in particular, since $B\in\Z$ we have $B\geq 4r$.
By a convexity argument one sees that $r\leq B\log 2$ implies $(B+1)^r\leq 2B^r$.
(Indeed, this follows immediately using the inequality $r\log(1+1/B)\leq r/B\leq \log 2$.)

Using all this, we have
$$
  \frac{1}{2}\left(\frac{4r}{B}\right)^r
  \leq
  \left(\frac{4r}{B+1}\right)^r
  \leq
  \frac{1}{p^{1/2}}
  \,,
$$
and hence
$$
  \sum_{x\in\F_p}
  \left|
  \sum_{1\leq b\leq B}
  \chi(x+b)
  \right|^{2r}
  \leq
  B^{2r}p^{1/2}\left(2r-\frac{1}{2}\right)
  \,.
$$

All together, this gives
$$
  \sum_{x\in\F_p}
  \nu_A(x)
  \left|
  \sum_{1\leq b\leq B}
  \chi(x+b)
  \right|
  \leq
    D_1\,
    \pi(A)^{1-\frac{1}{2r}}H^{1-\frac{1}{2r}}Bp^{\frac{1}{4r}}
  $$  
  with
  $$
      D_1=2^{\frac{1}{2r}}\left(2r-\frac{1}{2}\right)^{\frac{1}{2r}}
      =
      (4r-1)^{\frac{1}{2r}}
      \,.
  $$
Therefore
\begin{eqnarray*}
  \frac{1}{\pi(A)B}
  \sum_{x\in\F_p}
  \nu_A(x)
  \left|
  \sum_{1\leq b\leq B}
  \chi(x+b)
  \right|
  &\leq&
    D_1\,
    H^{1-\frac{1}{r}}p^{\frac{1}{4r}}
    \left(\frac{H}{\pi(A)}\right)^{\frac{1}{2r}}
    \,.
\end{eqnarray*}

Using (3.5) of~\cite{rosser.schoenfeld} and some simple computation,
provided $A\geq 3$ and $A\in\Z$, 
we have
$\pi(A)\geq A/(1+\log A)$;
using this, together with the bound
\begin{eqnarray*}
  \log A
  &\leq&
  \left(\frac{1}{2}-\frac{1}{4r}\right)\log p-\log (3\sqrt{rd})\\
  &<&
  \left(\frac{1}{2}-\frac{1}{4r}\right)\log p-1
  \,,
\end{eqnarray*}
allows us to estimate
\begin{eqnarray*}
  \frac{H}{\pi(A)}
  &\leq&
  \frac{H(\log A+1)}{A}
  \\
  &\leq&
  (d+1)B(\log A+1)
  \\
  &\leq&
  4r(d+1)\left(\frac{1}{2}-\frac{1}{4r}\right)p^{\frac{1}{2r}}\log p
  \,.
\end{eqnarray*}
%
%
Therefore
$$
  \left(  \frac{H}{\pi(A)}\right)^{\frac{1}{2r}}
  \leq
  D_2
  \,
  p^{\frac{1}{4r^2}}(\log p)^{\frac{1}{2r}}
$$
with
$$
  D_2=
  \left[
  4r(d+1)\left(\frac{1}{2}-\frac{1}{4r}\right)
  \right]^{\frac{1}{2r}}
  =
  (
  (d+1)(2r-1)
  )^{\frac{1}{2r}}
  \,,
$$
which leads to
\begin{eqnarray*}
  \frac{1}{\pi(A)B}
  \sum_{x\in\F_p}
  \nu_A(x)
  \left|
  \sum_{1\leq b\leq B}
  \chi(x+b)
  \right|
  &\leq&
    D_1 D_2\,
    H^{1-\frac{1}{r}}p^{\frac{r+1}{4r^2}}(\log p)^{\frac{1}{2r}}
    =
    D_1 D_2 E(H)
    \,.
\end{eqnarray*}
Finally,
using
(\ref{E:charsum.finally}),
this gives
$$
  |S_\chi(H)|
  \leq
  \left[
  \left(
    (d+1)(2r-1)(4r-1)
  \right)^{\frac{1}{2r}}
  +
  \frac{2C}{d^{1-\frac{1}{r}}}
  \right]
  E(H)
  \,.
$$
Now we see that $|S_\chi(H)|\leq C\,E(H)$, which would complete our induction, provided
\begin{equation}\label{E:equiv.hyp}
  \left(
  (d+1)(2r-1)(4r-1)
  \right)^{\frac{1}{2r}}
  +
  \frac{2C}{d^{1-\frac{1}{r}}}
  \leq
  C
  \,.
\end{equation}
Using the fact
$$
  d>4
  \quad\Longrightarrow\quad
  1-\frac{2}{d^{1-\frac{1}{r}}}>0
  \,,
$$
and solving (\ref{E:equiv.hyp}) for $C$
allows us to see that (\ref{E:equiv.hyp}) is equivalent to
(\ref{E:charsum.H2}).
\qed
%

\vspace{1ex}
\noindent\textbf{Proof of Theorem~\ref{T:burgess}.}
We apply Theorem~\ref{T:charsum} with $d=11$, $p_0=2\cdot 10^4$
and perform the necessary numerical computations, being careful to round up in our computations of values for $C(r)$.
\qed

\vspace{1ex}
The choices of $p_0$ and $d$ in the proof of Theorem~\ref{T:burgess} 
were designed to easily derive a widely applicable version of the
character sum estimate with decent constants for all $r$.
This will suit our purposes here.  However, if one wanted to achieve a slightly
better constant for a specific application, one would proceed as follows:
for any given $r$ and $p_0$, choose (or numerically estimate) the parameter $d$ so as to minimize~$C$.

\end{document}